\documentclass[11pt]{article}
\usepackage{geometry,graphicx,latexsym,amssymb,verbatim}
\usepackage{tikz}
\usetikzlibrary{calc}

\def\vertex(#1){\put(#1){\circle*{2}}}
\def\vertexo(#1){\put(#1){\circle{2}}}
\def\vert(#1){\put(#1){\circle*{1.5}}}
\def\verto(#1){\put(#1){\circle{1.5}}}
\def\lab(#1)#2{\put(#1){\makebox(0,0)[c]{#2}}}
\newtheorem{thm}{Theorem}
\newtheorem{ob}[thm]{Observation}
\newtheorem{lem}[thm]{Lemma}
\newtheorem{prop}[thm]{Proposition}
\newtheorem{cor}[thm]{Corollary}

\newcommand{\dtd}{\gamma^d_t}

\newcommand{\gt}{\gamma_t}
\newcommand{\diam}{{\rm diam}}
\newcommand{\mod}{{\rm mod}}
\newcommand{\qed}{$\Box$}

\newcommand{\cE}{{\cal E}}
\newcommand{\cF}{{\cal F}}
\newcommand{\cG}{{\cal G}}
\newcommand{\cH}{{\cal H}}
\newcommand{\cL}{{\cal L}}
\newcommand{\cT}{{\cal T}}
\newcommand{\cS}{{\cal S}}

\newcommand{\proof}{\noindent\textbf{Proof. }}
\newcommand{\smallqed}{{\tiny ($\Box$)}}

\newenvironment{unnumbered}[1]{\trivlist
\item [\hskip \labelsep {\bf #1}]\ignorespaces\it}{\endtrivlist}

\newcommand{\2}{\vspace{0.2cm}}

\begin{document}

\title{Disjunctive Total Domination in Graphs}
\author{$^1$Michael A. Henning\thanks{Research supported in part
by the South African National Research Foundation and the University of Johannesburg}  \, and $^{1,2}$Viroshan Naicker \\
\\
$^1$Department of Mathematics\\
University of Johannesburg \\
Auckland Park, 2006 South Africa \\
Email: mahenning@uj.ac.za\\
\\
$^2$Department of Mathematics \\
Rhodes University \\
Grahamstown, 6140 South Africa \\
Email: V.Naicker@ru.ac.za
 }

\date{}
\maketitle

\begin{abstract}
Let $G$ be a graph with no isolated vertex. In this paper, we study a parameter that is a relaxation of arguably the most important domination parameter, namely the total domination number, $\gamma_t(G)$. A set $S$ of vertices in $G$ is a disjunctive total dominating set of $G$ if every vertex is adjacent to a vertex of $S$ or has at least two vertices in $S$ at distance~$2$ from it. The disjunctive total domination number, $\dtd(G)$, is the minimum cardinality of such a set. We observe that $\dtd(G) \le \gamma_t(G)$. We prove that if $G$ is a connected graph of order~$n \ge 8$, then $\dtd(G) \le 2(n-1)/3$ and we characterize the extremal graphs. It is known that if $G$ is a connected claw-free graph of order~$n$, then $\gamma_t(G) \le 2n/3$ and this upper bound is tight for arbitrarily large~$n$. We show this upper bound can be improved significantly for the disjunctive total domination number. We show that if $G$ is a connected claw-free graph of order~$n > 10$, then $\dtd(G) \le 4n/7$ and we characterize the graphs achieving equality in this bound.
\end{abstract}

{\small \textbf{Keywords:} Total dominating set; disjunctive total dominating set; claw-free. }\\
\indent {\small \textbf{AMS subject classification: 05C69}}

\newpage
\section{Introduction}

A major issue in modern networks is to minimize the trade-off between resource allocation and redundancy. Key resources are usually expensive and cannot be allocated across an entire network, and, in addition, if there is a possibility of resource failure at a particular node, redundancy and backup requirements then become vital but require extra resources to be allocated. This problem has been addressed, in various guises, by using graphs as a model for the network and searching for vertex subsets which are `close' to the rest of the graph and satisfy pertinent redundancy criteria.

Domination and, in particular, total domination are well studied topics in the graph theory literature which attempt a solution of this problem~(see, for example, \cite{hhs1,hhs2,He09,HeYe_book}).

Suppose that $G$ is a graph with vertex $V$ that serves as a model of a network. On the one hand, for purposes of resource allocation, we select a set $D$ of vertices, called a \emph{dominating set}, of $G$ such that every vertex in $V \setminus D$ is adjacent to at least one vertex in $D$. On the other hand, for the purpose of extending the domination problem to include redundancy, we select a set $S$ of vertices, called a \emph{total dominating set} and abbreviated TD-set, of $G$ such that such that every vertex in $V$, including those in $S$, is adjacent to at least one vertex in $S$.
However, given the sheer scale of modern networks, many existing domination type structures are expensive to implement. Variations on the theme of dominating and total dominating sets studied to date tend to focus on adding restrictions which in turn raises their implementation costs.  As an alternative route a relaxation of the domination number, called \emph{disjunctive domination}, was proposed and studied by Goddard et al.~\cite{GoHePi11}. In this paper we extend this concept to a relaxation of total domination, called \emph{disjunctive total domination}, which allows for greater flexibility in the modeling networks where one trades off redundancy and backup capability with resource optimization.

\begin{figure}[htb]
\tikzstyle{every node}=[circle, draw, fill=black!0, inner sep=0pt,minimum width=.18cm]
\begin{center}
\begin{tikzpicture}[thick,scale=.6]
  \draw(0,0) { 
    +(0.91,1.82) -- +(1.82,1.82)
    +(1.82,1.82) -- +(2.73,0.91)
    +(2.73,0.91) -- +(1.82,0.00)
    +(1.82,0.00) -- +(0.91,0.00)
    +(0.91,0.00) -- +(0.00,0.91)
    +(0.00,0.91) -- +(0.91,1.82)
    +(2.73,0.91) -- +(3.64,0.91)
    +(3.64,0.91) -- +(4.55,1.82)
    +(4.55,1.82) -- +(5.45,0.91)
    +(5.45,0.91) -- +(4.55,0.00)
    +(4.55,0.00) -- +(3.64,0.91)
    +(7.27,0.91) -- +(8.18,1.82)
    +(8.18,1.82) -- +(9.09,0.91)
    +(8.18,0.00) -- +(9.09,0.91)
    +(7.27,0.91) -- +(8.18,0.00)
    +(5.45,0.91) -- +(6.36,0.91)
    +(6.36,0.91) -- +(7.27,0.91)
    +(0.91,1.82) node{}
    +(1.82,1.82) node[circle, draw,fill=black!100]{}
    +(1.82,0.00) node[circle, draw,fill=black!100]{}
    +(0.91,0.00) node{}
    +(0.00,0.91) node{}
    +(2.73,0.91) node{}
    +(3.64,0.91) node[circle, draw,fill=black!100]{}
    +(4.55,1.82) node{}
    +(4.55,0.00) node{}
    +(5.45,0.91) node{}
    +(7.27,0.91) node[circle, draw,fill=black!100]{}
    +(8.18,1.82) node[circle, draw,fill=black!100]{}
    +(8.18,0.00) node{}
    +(9.09,0.91) node{}
    +(6.36,0.91) node{}
  };
\end{tikzpicture}
\end{center}
\vskip -0.6 cm \caption{A graph $G$ with $\dtd(G) = 5$.} \label{f:DTD5}
\end{figure}
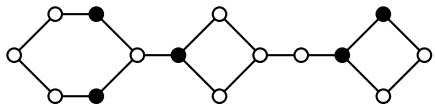

A set $S$ of vertices in $G$ is a \emph{disjunctive total dominating set}, abbreviated DTD-set, of $G$ if every vertex is adjacent to a vertex of $S$ or has at least two vertices in $S$ at distance~$2$ from it. For example, the set of five darkened vertices in the graph $G$ shown in Figure~\ref{f:DTD5} is a DTD-set of $G$. We will say that a vertex $v$ is \emph{disjunctively totally dominated} by a set $S$ if $v$ has a neighbor in $S$ or if $v$ is at distance~$2$ from at least two vertices of $S$. The \emph{disjunctive total domination number}, $\dtd(G)$, is the minimum cardinality of a DTD-set in $G$. A DTD-set of cardinality $\dtd(G)$ is called a $\dtd(G)$-set.

The \emph{domination number} of $G$, denoted $\gamma(G)$, is the minimum cardinality of a dominating set in $G$, while the \emph{total domination number} of $G$, denoted by $\gt(G)$, is the minimum cardinality of a TD-set of $G$. Every TD-set is a DTD-set, implying the following observation.

\begin{ob}
For every graph $G$ with no isolated vertex, $\dtd(G) \le \gt(G)$.
 \label{ob:Tdom}
\end{ob}

\newpage
\subsection{Notation}

For notation and graph theory terminology, we in general
follow~\cite{hhs1}. Specifically, let $G = (V,E)$ be a graph with
vertex set $V$, edge set $E$ and no isolated vertex. The \emph{open
neighborhood} of a vertex $v \in V$ is $N_G(v) = \{u \in V \, | \, uv
\in E(G)\}$ and its \emph{closed neighborhood} is the set $N_G[v] =
N_G(v) \cup \{v\}$. The degree of $v$ is $d_G(v) = |N_G(v)|$.
The minimum degree among the vertices of $G$ is denoted by $\delta(G)$.
For a set $S \subseteq V$, the subgraph induced by $S$ is denoted by $G[S]$, while the graph obtained from $G$ be removing all vertices in $S$ and their incident edges is denoted by $G - S$. For two vertices $u$ and $v$ in a connected graph $G$, the \emph{distance} $d_G(u,v)$ between $u$ and $v$ is the length of a shortest $u$--$v$ path in $G$.
If the graph $G$ is clear from the context, we simply write $N(v)$, $N[v]$ and $d(v)$ rather than $N_G(v)$, $N_G[v]$  and $d_G(v)$,
respectively.

For a set $S \subseteq V$, its \emph{open neighborhood}
is the set $N(S) = \cup_{v \in S} N(v)$ and its \emph{closed
neighborhood} is the set $N[S] = N(S) \cup S$. For sets $A, B \subseteq V$, we say that $A$ \emph{totally dominates} $B$ if $B \subseteq N(A)$. A \emph{cycle} and \emph{path} on $n$ vertices are denoted by $C_n$ and $P_n$, respectively.

For a graph $H$ and integer $k \ge 1$, we denote by $H \circ P_k$ the
graph of order~$(k+1)|V(H)|$ obtained from $H$ by attaching a path of
length~$k$ to each vertex of $H$ so that the resulting paths are
vertex-disjoint. The graph $H \circ P_2$ is also called the
\emph{$2$-corona} of $H$.

A \emph{leaf} is a vertex of degree~$1$ and a \emph{support vertex} is a vertex adjacent to a leaf. A \emph{star}
is a tree $K_{1,n}$, while for integers $r,s \ge 1$, a \emph{double
star} $S(r,s)$ is a tree with exactly two vertices that are not
leaves, one of which has degree~$r+1$ and the other degree~$s+1$.
A \emph{rooted tree} distinguishes one vertex $r$ called the
\emph{root}. For each vertex $v \ne r$ of $T$, the \emph{parent} of
$v$ is the neighbor of $v$ on the unique $r$--$v$ path, while a
\emph{child} of $v$ is any other neighbor of~$v$.

\subsection{Known Results}

Cockayne, Dawes, and Hedetniemi~\cite{CoDaHe80} proved that the total domination number of a connected graph of order at least~$3$ is bounded above by two-thirds its order. The graphs achieving equality in this bound were characterized by Brigham, Carrington, and Vitray~\cite{BrCaVi00}. If we restrict $G$ to be a connected claw-free graph, then this upper bound of two-thirds the order of the graph cannot be improved since the $2$-corona of a complete graph is claw-free and has total domination number two-thirds its order.

\begin{thm}{\rm (\cite{BrCaVi00,CoDaHe80})}
If $G$ is a connected graph of order~$n \ge 3$, then $\gt(G) \le
2n/3$. Further, this bound is sharp even if we restrict our attention to the class of claw-free graphs.
 \label{claw_free_known}
\end{thm}

Let $C_{10}'$ be the graph obtained from a $10$-cycle $v_1v_2 \ldots v_{10}v_1$ by adding the edge $v_1v_6$ and let $C_{10}''$ be the graph obtained from $C_{10}'$ by adding the edge $v_2v_7$. If we restrict the minimum degree to be two and impose the additional restriction that the graph is connected, the upper bound on the total domination number can be improved as follows.

\begin{thm}{\rm (\cite{He00})}
If $G \not\in \{C_3, C_5, C_6, C_{10}, C_{10}', C_{10}''\}$ is a
connected graph of order~$n$ with $\delta(G) \ge 2$, then $\gt(G) \le
4n/7$.
 \label{Tdom_deg2}
\end{thm}

\begin{figure}[htb]
\tikzstyle{every node}=[circle, draw, fill=black!0, inner sep=0pt,minimum width=.16cm]
\begin{center}
\begin{tikzpicture}[thick,scale=.7]

  \draw(0,15) {
    +(0.00,0.00)--+(0.50,0.00) +(3.00,0.00)--+(2.50,0.00) +(1.00,0.00)--+(0.50,0.00) +(1.00,0.00)--+(1.50,0.00) +(2.00,0.00)--+(2.50,0.00) +(2.00,0.00)--+(1.50,0.00)
    +(0.50,0.00) node{} +(2.50,0.00) node{} +(1.50,0.00) node{} +(0.00,0.00) node{} +(3.00,0.00) node{} +(1.00,0.00) node{} +(2.00,0.00) node{}
  +(1.5,-0.75) node[rectangle, draw=white!0, fill=white!100]{(a) $L_{1}$}
  };
  \draw(5,15) {
    +(1.25,0.50)--+(1.50,0.00) +(0.00,0.00)--+(0.50,0.00) +(2.00,0.00)--+(1.50,0.00) +(2.00,0.00)--+(2.50,0.00) +(1.00,0.00)--+(1.50,0.00) +(1.00,0.00)--+(0.50,0.00) +(1.00,0.00)--+(1.25,0.50)
    +(1.50,0.00) node{} +(0.50,0.00) node{} +(2.50,0.00) node{} +(1.25,0.50) node{} +(0.00,0.00) node{} +(2.00,0.00) node{} +(1.00,0.00) node{}
  +(1.5,-0.75) node[rectangle, draw=white!0, fill=white!100]{(b) $L_{2}$}
  };
  \draw(10,15) {
    +(2.00,0.00)--+(1.75,0.50) +(2.00,0.00)--+(2.50,0.00) +(0.50,0.00)--+(1.00,0.00) +(0.50,0.00)--+(0.00,0.00) +(1.50,0.00)--+(1.75,0.50) +(1.50,0.00)--+(1.00,0.00) +(1.50,0.00)--+(2.00,0.00)
    +(1.75,0.50) node{} +(2.50,0.00) node{} +(1.00,0.00) node{} +(0.00,0.00) node{} +(2.00,0.00) node{} +(0.50,0.00) node{} +(1.50,0.00) node{}
  +(1.5,-0.75) node[rectangle, draw=white!0, fill=white!100]{(c) $L_{3}$}
  };
  \draw(15,15) {
    +(1.25,1.00)--+(1.25,0.50) +(0.00,0.00)--+(0.50,0.00) +(1.50,0.00)--+(1.25,0.50) +(1.50,0.00)--+(2.00,0.00) +(1.00,0.00)--+(1.25,0.50) +(1.00,0.00)--+(0.50,0.00) +(1.00,0.00)--+(1.50,0.00)
    +(1.25,0.50) node{} +(0.50,0.00) node{} +(2.00,0.00) node{} +(1.25,1.00) node{} +(0.00,0.00) node{} +(1.50,0.00) node{} +(1.00,0.00) node{}
  +(1,-0.75) node[rectangle, draw=white!0, fill=white!100]{(d) $L_{4}$}
  };
    \draw(0,12) {
    +(0.00,0.50)--+(0.50,0.50) +(2.50,0.00)--+(2.50,1.00) +(1.00,0.50)--+(0.50,0.50) +(1.00,0.50)--+(1.50,0.50) +(2.00,0.50)--+(2.50,1.00) +(2.00,0.50)--+(1.50,0.50) +(2.00,0.50)--+(2.50,0.00)
    +(0.50,0.50) node{} +(2.50,1.00) node{} +(1.50,0.50) node{} +(0.00,0.50) node{} +(2.50,0.00) node{} +(1.00,0.50) node{} +(2.00,0.50) node{}
  +(1.5,-0.75) node[rectangle, draw=white!0, fill=white!100]{(e) $L_{5}$}
   };
   \draw(5,12) {
    +(2.00,1.00)--+(2.00,0.00) +(2.50,0.50)--+(2.00,0.00) +(2.50,0.50)--+(2.00,1.00) +(0.50,0.50)--+(1.00,0.50) +(0.50,0.50)--+(0.00,0.50) +(1.50,0.50)--+(2.00,0.00) +(1.50,0.50)--+(1.00,0.50) +(1.50,0.50)--+(2.00,1.00)
    +(2.00,0.00) node{} +(1.00,0.50) node{} +(0.00,0.50) node{} +(2.00,1.00) node{} +(2.50,0.50) node{} +(0.50,0.50) node{} +(1.50,0.50) node{}
  +(1,-0.75) node[rectangle, draw=white!0, fill=white!100]{(f) $L_{6}$}
  };
    \draw(10,12) {
    +(0.00,0.50)--+(0.50,0.50) +(1.00,0.50)--+(0.50,0.50) +(1.00,0.50)--+(1.50,0.00) +(2.50,0.00)--+(1.50,0.00) +(2.50,0.00)--+(2.50,1.00) +(1.50,1.00)--+(1.50,0.00) +(1.50,1.00)--+(2.50,1.00) +(1.50,1.00)--+(1.00,0.50)
    +(0.50,0.50) node{} +(1.50,0.00) node{} +(2.50,1.00) node{} +(0.00,0.50) node{} +(1.00,0.50) node{} +(2.50,0.00) node{} +(1.50,1.00) node{}
  +(1.5,-0.75) node[rectangle, draw=white!0, fill=white!100]{(g) $L_{7}$} };
   \draw(15,12) {
    +(0.00,0.50)--+(0.50,0.50) +(1.00,0.00)--+(0.50,0.50) +(1.00,0.00)--+(1.50,0.50) +(2.00,0.50)--+(1.50,0.50) +(2.00,0.50)--+(2.50,0.50) +(1.00,1.00)--+(0.50,0.50) +(1.00,1.00)--+(1.50,0.50) +(1.00,1.00)--+(1.00,0.00)
    +(0.50,0.50) node{} +(1.50,0.50) node{} +(2.50,0.50) node{} +(0.00,0.50) node{} +(1.00,0.00) node{} +(2.00,0.50) node{} +(1.00,1.00) node{}
  +(1.5,-0.75) node[rectangle, draw=white!0, fill=white!100]{(h) $L_{8}$}
 };
   \draw(0.5,7.5) {
    +(0.95,3.00)--+(0.95,2.50) +(1.54,1.81)--+(0.95,2.50) +(1.54,1.81)--+(1.90,0.69) +(0.95,0.69)--+(1.90,0.69) +(0.95,0.69)--+(0.00,0.69) +(0.36,1.81)--+(0.95,2.50) +(0.36,1.81)--+(0.00,0.69) +(0.36,1.81)--+(1.54,1.81)
    +(0.95,2.50) node{} +(1.90,0.69) node{} +(0.00,0.69) node{} +(0.95,3.00) node{} +(1.54,1.81) node{} +(0.95,0.69) node{} +(0.36,1.81) node{}
  +(1,0) node[rectangle, draw=white!0, fill=white!100]{(i) $L_{9}$}
 };
  \draw(5.5,8) {
    +(0.00,0.68)--+(0.19,1.52) +(0.97,1.90)--+(0.19,1.52) +(0.97,1.90)--+(1.76,1.52) +(1.95,0.68)--+(1.76,1.52) +(1.95,0.68)--+(1.41,0.00) +(0.54,0.00)--+(1.41,0.00) +(0.54,0.00)--+(0.00,0.68)
    +(0.19,1.52) node{} +(1.76,1.52) node{} +(1.41,0.00) node{} +(0.00,0.68) node{} +(0.97,1.90) node{} +(1.95,0.68) node{} +(0.54,0.00) node{}
  +(1,-0.75) node[rectangle, draw=white!0, fill=white!100]{(j) $L_{10}$}
  };
  \draw(10.5,8) {
    +(0.87,1.00)--+(0.00,1.50) +(0.87,1.00)--+(1.73,1.50) +(0.00,0.50)--+(0.00,1.50) +(0.00,0.50)--+(0.87,0.00) +(1.73,0.50)--+(1.73,1.50) +(1.73,0.50)--+(0.87,0.00) +(0.87,2.00)--+(0.00,1.50) +(0.87,2.00)--+(1.73,1.50) +(0.87,2.00)--+(0.87,1.00)
    +(0.00,1.50) node{} +(1.73,1.50) node{} +(0.87,0.00) node{} +(0.87,1.00) node{} +(0.00,0.50) node{} +(1.73,0.50) node{} +(0.87,2.00) node{}
  +(1,-0.75) node[rectangle, draw=white!0, fill=white!100]{(k) $L_{11}$}
  };
  \draw(15,8) {
    +(1.00,2.37)--+(1.50,1.73) +(2.00,0.87)--+(1.50,1.73) +(2.00,0.87)--+(1.50,0.00) +(0.50,0.00)--+(1.50,0.00) +(0.50,0.00)--+(0.00,0.87) +(0.50,1.73)--+(1.50,1.73) +(0.50,1.73)--+(0.00,0.87) +(0.50,1.73)--+(1.00,2.37)
    +(1.50,1.73) node{} +(1.50,0.00) node{} +(0.00,0.87) node{} +(1.00,2.37) node{} +(2.00,0.87) node{} +(0.50,0.00) node{} +(0.50,1.73) node{}
    +(1,-0.75) node[rectangle, draw=white!0, fill=white!100]{($\ell$) $L_{12}$}
  };
\end{tikzpicture}
\end{center}
\vskip -0.5 cm \caption{The $12$ claw-free graphs on $7$ vertices
with total domination number $4$.} \label{f:cL}
\end{figure}
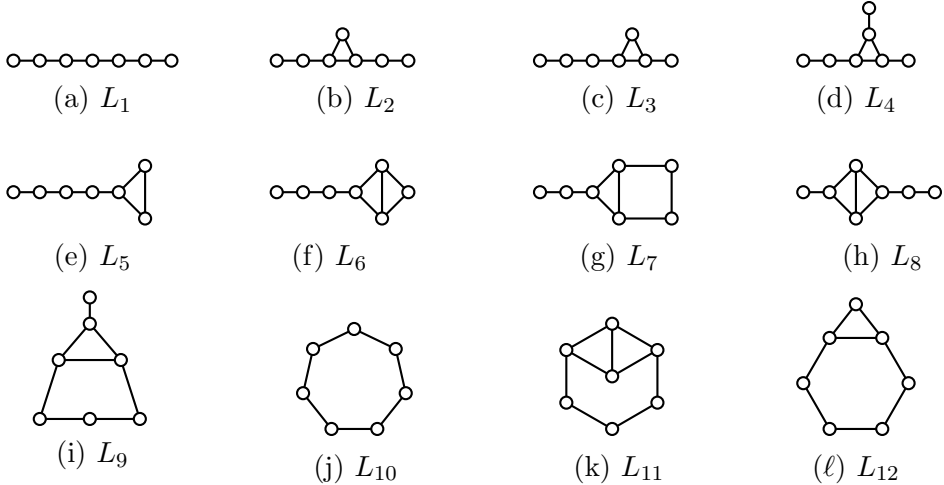

The total domination number of a connected graph on seven vertices is
at most~$4$. It is a simple exercise to check that there are exactly
20 such graphs with total domination number~$4$. Of these 20 graphs,
12 are claw-free and belong to the family $\cL =
\{L_1,L_2,\ldots,L_{12}\}$ shown in Figure~\ref{f:cL}. Hence we have
the following observation. Note that $L_1 = P_7$ and $L_{10} = C_7$.

\begin{ob}
Let $G$ be a connected claw-free graph of order~$n = 7$. Then,
$\gt(G) \le 4$ with equality if and only if $G \in \cL$.
 \label{ob:n7}
\end{ob}

A graph $G$ is referred to in~\cite{He00} as a
$\frac{4}{7}$-\emph{minimal graph} if $G$ is edge-minimal with
respect to satisfying the following three conditions: (i) $\delta(G)
\ge 2$, (ii) $G$ is connected, and (iii) $\gt(G) \ge 4n/7$, where $n$
is the order of $G$. The $\frac{4}{7}$-minimal graphs are
characterized in~\cite{He00}. As a consequence of this
characterization and of the main result in~\cite{He00} which
characterizes the connected graphs $G$ of order~$n > 14$  with
$\delta(G) \ge 2$ achieving the upper bound of
Theorem~\ref{Tdom_deg2}, one can readily obtain the following result.

\begin{cor}
Let $G$ be a connected claw-free graph of order~$n$ with $\delta(G)
\ge 2$. If $G \notin \{C_3, C_5, C_6, C_{10} \}$, then $\gt(G) \le
4n/7$. Further if $\gt(G) = 4n/7$, then $G \in \{L_{10}, L_{11},
L_{12}, C_{14}\}$. \label{ob:clawfree}
\end{cor}

\subsection{Special Families}
\label{S:special}

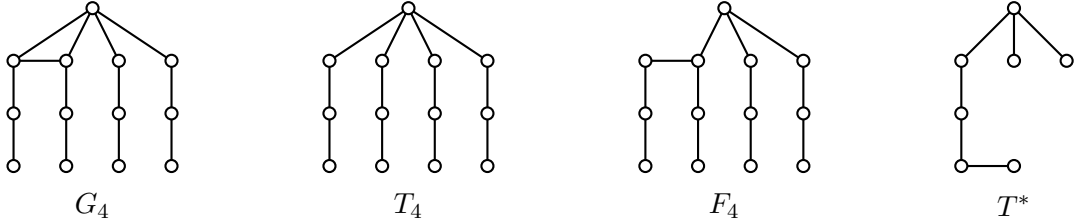
\begin{figure}[htb]
\tikzstyle{every node}=[circle, draw, fill=black!0, inner sep=0pt,minimum width=.16cm]
\begin{center}
\begin{tikzpicture}[thick,scale=.7]
\draw(0,0) { 
    +(1.50,2.00) -- +(0.00,1.00)
    +(0.00,1.00) -- +(1.00,1.00)
    +(0.00,1.00) -- +(0.00,0.00)
    +(0.00,0.00) -- +(0.00,-1.00)
    +(1.00,0.00) -- +(1.00,-1.00)
    +(2.00,0.00) -- +(2.00,-1.00)
    +(3.00,0.00) -- +(3.00,-1.00)
    +(1.00,0.00) -- +(1.00,1.00)
    +(1.00,1.00) -- +(1.50,2.00)
    +(1.50,2.00) -- +(2.00,1.00)
    +(2.00,1.00) -- +(2.00,0.00)
    +(3.00,0.00) -- +(3.00,1.00)
    +(3.00,1.00) -- +(1.50,2.00)
    +(0.00,-1.00) node{}
    +(1.00,-1.00) node{}
    +(2.00,-1.00) node{}
    +(3.00,-1.00) node{}
    +(0.00,1.00) node{}
    +(1.00,1.00) node{}
    +(2.00,1.00) node{}
    +(3.00,1.00) node{}
    +(1.50,2.00) node{}
    +(0.00,0.00) node{}
    +(1.00,0.00) node{}
    +(2.00,0.00) node{}
    +(3.00,0.00) node{}
    +(1.50,-1.75) node[rectangle, draw=white!0, fill=white!100]{$G_4$}
   };
  \draw(6,0) { 
    +(1.50,2.00) -- +(0.00,1.00)
    +(0.00,1.00) -- +(0.00,0.00)
    +(1.00,0.00) -- +(1.00,1.00)
    +(1.00,1.00) -- +(1.50,2.00)
    +(1.50,2.00) -- +(2.00,1.00)
    +(2.00,1.00) -- +(2.00,0.00)
    +(3.00,0.00) -- +(3.00,1.00)
    +(3.00,1.00) -- +(1.50,2.00)
    +(6.00,1.00) -- +(6.00,0.00)
    +(7.00,0.00) -- +(7.00,1.00)
    +(7.00,1.00) -- +(7.50,2.00)
    +(7.50,2.00) -- +(8.00,1.00)
    +(8.00,1.00) -- +(8.00,0.00)
    +(9.00,0.00) -- +(9.00,1.00)
    +(6.00,1.00) -- +(7.00,1.00)
    +(9.00,1.00) -- +(7.50,2.00)
    +(0.00,0.00) -- +(0.00,-1.00)
    +(1.00,0.00) -- +(1.00,-1.00)
    +(2.00,0.00) -- +(2.00,-1.00)
    +(3.00,0.00) -- +(3.00,-1.00)
    +(6.00,0.00) -- +(6.00,-1.00)
    +(7.00,0.00) -- +(7.00,-1.00)
    +(8.00,0.00) -- +(8.00,-1.00)
    +(9.00,0.00) -- +(9.00,-1.00)
    +(0.00,-1.00) node{}
    +(1.00,-1.00) node{}
    +(2.00,-1.00) node{}
    +(3.00,-1.00) node{}
    +(0.00,1.00) node{}
    +(1.00,1.00) node{}
    +(2.00,1.00) node{}
    +(3.00,1.00) node{}
    +(1.50,2.00) node{}
    +(0.00,0.00) node{}
    +(1.00,0.00) node{}
    +(2.00,0.00) node{}
    +(3.00,0.00) node{}
    +(6.00,1.00) node{}
    +(7.00,1.00) node{}
    +(8.00,1.00) node{}
    +(9.00,1.00) node{}
    +(6.00,0.00) node{}
    +(7.00,0.00) node{}
    +(8.00,0.00) node{}
    +(9.00,0.00) node{}
    +(7.50,2.00) node{}
    +(6.00,-1.00) node{}
    +(7.00,-1.00) node{}
    +(8.00,-1.00) node{}
    +(9.00,-1.00) node{}
    +(1.50,-1.75) node[rectangle, draw=white!0, fill=white!100]{$T_4$}
    +(7.50,-1.75) node[rectangle, draw=white!0, fill=white!100]{$F_4$}
  };
  \draw(18,-1) { 
    +(1.00,3.00) -- +(1.00,2.00)
    +(1.00,3.00) -- +(2.00,2.00)
    +(1.00,3.00) -- +(0.00,2.00)
    +(0.00,2.00) -- +(0.00,1.00)
    +(0.00,1.00) -- +(0.00,0.00)
    +(0.00,0.00) -- +(1.00,0.00)
    +(0.00,1.00) node{}
    +(0.00,0.00) node{}
    +(1.00,0.00) node{}
    +(0.00,2.00) node{}
    +(1.00,3.00) node{}
    +(2.00,2.00) node{}
    +(1.00,2.00) node{}
    +(1.00,-0.75) node[rectangle, draw=white!0, fill=white!100]{$T^*$}   
  };
\end{tikzpicture}
\end{center}
\vskip -0.6 cm \caption{The graph $G_4$ and the trees $T_4$, $F_4$ and $T^*$.} \label{f:tree}
\end{figure}

For $k \ge 1$, let $T_k$ be the tree obtained from a star $K_{1,k}$
by subdividing every edge exactly twice and let $\cT$ be the family
of all such trees $T_k$. For $k \ge 2$, let $G_k$ be the graph
obtained from $T_k$ by adding an edge joining two neighbors of the
central vertex of $T_k$, and let $\cG$ be the family of all such
graphs $G_k$. Finally for $k \ge 2$, let $F_k$ be the tree obtained
from $T_k$ by deleting an edge $uv$ incident with the central vertex
$v$ of $T_k$ and adding the edge $uw$ for some neighbor $w$ of $v$
different from $u$, and let $\cF$ be the family of all such trees
$F_k$. Let $T^*$ be the tree obtained from a star $K_{1,3}$ by
subdividing one edge three times, and so $T^*$ has order~$7$. The
graph $G_4$ and the trees $T_4$, $F_4$ and $T^*$ are illustrated in
Figure~\ref{f:tree}.

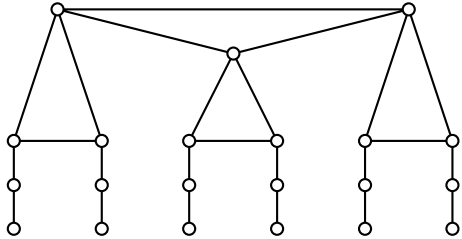
\begin{figure}[htb]
\tikzstyle{every node}=[circle, draw, fill=black!0, inner sep=0pt,minimum width=.16cm]
\begin{center}
\begin{tikzpicture}[thick,scale=.7]
  \draw(0,0) { 
    +(0.83,4.17) -- +(0.00,1.67)
    +(0.00,1.67) -- +(1.67,1.67)
    +(1.67,1.67) -- +(0.83,4.17)
    +(0.83,4.17) -- +(7.50,4.17)
    +(0.83,4.17) -- +(4.17,3.33)
    +(4.17,3.33) -- +(7.50,4.17)
    +(4.17,3.33) -- +(3.33,1.67)
    +(3.33,1.67) -- +(5.00,1.67)
    +(5.00,1.67) -- +(4.17,3.33)
    +(7.50,4.17) -- +(6.67,1.67)
    +(6.67,1.67) -- +(8.33,1.67)
    +(8.33,1.67) -- +(7.50,4.17)
    +(0.00,1.67) -- +(0.00,0.83)
    +(0.00,0.83) -- +(0.00,0.00)
    +(1.67,1.67) -- +(1.67,0.83)
    +(1.67,0.83) -- +(1.67,0.00)
    +(3.33,0.00) -- +(3.33,0.83)
    +(3.33,0.83) -- +(3.33,1.67)
    +(5.00,1.67) -- +(5.00,0.83)
    +(5.00,0.83) -- +(5.00,0.00)
    +(6.67,0.00) -- +(6.67,0.83)
    +(6.67,0.83) -- +(6.67,1.67)
    +(8.33,1.67) -- +(8.33,0.83)
    +(8.33,0.83) -- +(8.33,0.00)
    +(0.83,4.17) node{}
    +(0.00,1.67) node{}
    +(0.00,0.83) node{}
    +(0.00,0.00) node{}
    +(1.67,1.67) node{}
    +(1.67,0.83) node{}
    +(1.67,0.00) node{}
    +(3.33,1.67) node{}
    +(3.33,0.83) node{}
    +(3.33,0.00) node{}
    +(5.00,0.00) node{}
    +(5.00,1.67) node{}
    +(5.00,0.83) node{}
    +(7.50,4.17) node{}
    +(6.67,1.67) node{}
    +(6.67,0.83) node{}
    +(6.67,0.00) node{}
    +(8.33,0.00) node{}
    +(8.33,0.83) node{}
    +(8.33,1.67) node{}
    +(4.17,3.33) node{}
  };
\end{tikzpicture}
\end{center}
\vskip -0.6 cm \caption{The graph $H_3$.} \label{f:H3}
\end{figure}

Let $H_t$ be the claw-free graph obtained from a complete graph $K_t$
on $t \ge 1$ vertices as follows: For each vertex $v$ of the complete
graph $K_t$, add a path $P_6$ and join $v$ to the two central
vertices of the path. Let $\cH$ be the family of all such graphs
$H_t$, where $t \ge 1$. The graph $H_3 \in \cH$ is illustrated in
Figure~\ref{f:H3}.

\medskip

\section{Results}

Our aim in this paper is twofold. First to establish a tight upper
bound on the disjunctive total domination number of a graph in terms
of its order and to characterize the extremal graphs. Secondly to
show that this bound can be significantly improved if we impose the
condition of claw-freeness on a graph. The following result
establishes an upper bound on the disjunctive total
domination number of a graph of order at least~$8$ and characterizes the graphs achieving equality in this bound. A proof of Theorem~\ref{char:graph} is given in Section~\ref{S:char:graph}.

\begin{thm}
Let $G$ be a connected graph of order~$n \ge 8$. Then, $\dtd(G) \le
2(n-1)/3$, with equality if and only if $G \in \cT \cup \cF \cup
\cG$.
 \label{char:graph}
\end{thm}

\begin{figure}[htb]
\tikzstyle{every node}=[circle, draw, fill=black!0, inner sep=0pt,minimum width=.16cm]
\begin{center}
\begin{tikzpicture}[thick,scale=.7]
  \draw(0,0) { 
    +(6.67,1.11) -- +(6.67,0.00)
    +(6.67,0.00) -- +(7.78,0.00)
    +(7.78,0.00) -- +(7.78,1.11)
    +(7.78,1.11) -- +(6.67,1.11)
    +(6.67,1.11) -- +(7.78,0.00)
    +(7.78,1.11) -- +(6.67,0.00)
    +(6.67,1.11) -- +(6.11,1.11)
    +(6.11,1.11) -- +(5.56,1.11)
    +(5.00,0.00) -- +(5.56,0.00)
    +(5.56,0.00) -- +(6.11,0.00)
    +(6.11,0.00) -- +(6.67,0.00)
    +(7.78,0.00) -- +(8.33,0.00)
    +(8.33,0.00) -- +(8.89,0.00)
    +(8.89,0.00) -- +(9.44,0.00)
    +(8.89,1.11) -- +(8.33,1.11)
    +(8.33,1.11) -- +(7.78,1.11)
    +(0.56,1.11) -- +(1.11,1.11)
    +(1.11,1.11) -- +(1.67,1.11)
    +(1.67,1.11) -- +(1.67,0.00)
    +(1.67,0.00) -- +(2.22,0.56)
    +(2.22,0.56) -- +(1.67,1.11)
    +(2.78,1.11) -- +(2.22,0.56)
    +(2.22,0.56) -- +(2.78,0.00)
    +(2.78,0.00) -- +(2.78,1.11)
    +(2.78,1.11) -- +(3.33,1.11)
    +(3.33,1.11) -- +(3.89,1.11)
    +(3.89,0.00) -- +(3.33,0.00)
    +(3.33,0.00) -- +(2.78,0.00)
    +(1.67,0.00) -- +(1.11,0.00)
    +(1.11,0.00) -- +(0.56,0.00)
    +(0.56,0.00) -- +(0.00,0.00)
    +(0.00,0.00) node{}
    +(0.56,0.00) node{}
    +(1.11,0.00) node{}
    +(1.67,0.00) node{}
    +(1.67,1.11) node{}
    +(2.22,0.56) node{}
    +(1.11,1.11) node{}
    +(0.56,1.11) node{}
    +(2.78,0.00) node{}
    +(3.33,0.00) node{}
    +(3.89,0.00) node{}
    +(2.78,1.11) node{}
    +(3.33,1.11) node{}
    +(3.89,1.11) node{}
    +(9.44,0.00) node{}
    +(8.89,0.00) node{}
    +(8.33,0.00) node{}
    +(7.78,0.00) node{}
    +(7.78,1.11) node{}
    +(8.33,1.11) node{}
    +(8.89,1.11) node{}
    +(6.67,1.11) node{}
    +(6.67,0.00) node{}
    +(6.11,1.11) node{}
    +(5.56,1.11) node{}
    +(6.11,0.00) node{}
    +(5.56,0.00) node{}
    +(5.00,0.00) node{}
    +(2,-0.75) node[rectangle, draw=white!0, fill=white!100]{(a) $L_{13}$}
    +(7.2,-0.75) node[rectangle, draw=white!0, fill=white!100]{(b) $L_{14}$}
  };
\end{tikzpicture}
\end{center}
\vskip -0.6 cm \caption{The graphs $L_{13}$ and $L_{14}$.} \label{f:L13}
\end{figure}
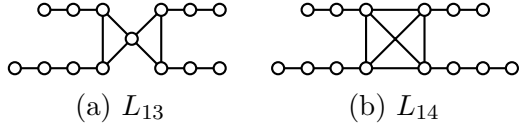

Total domination in claw-free graphs is well studied in the literature, see for example~\cite{FaHe03,FaHe08,FaHe08b,HeYe_claw_free,SoHe10,Li11} and elsewhere. We show next that if we impose the condition of claw-freeness on a
graph, then the upper bound in Theorem~\ref{claw_free_known} can be improved significantly from approximately two-thirds its order to four-sevenths its order. For this purpose, let $\cE = \{ P_2,P_3,P_5,P_6,C_3,G_3 \}$. We call a graph in the family $\cE$ an ``exceptional graph".
%
Let $\cS_1 = \{L_1,L_2,L_3,L_5,L_6,L_{10}\}$ and let
$\cS = \cS_1 \cup \{L_{13},L_{14}\}$, where
$L_{13}$ and $L_{14}$ are the graphs shown in Figure~\ref{f:L13}.

\begin{thm}
Let $G$ be a connected claw-free graph of order~$n \ge 2$. Then,
$\dtd(G) \le 4n/7$, unless $G \in \cE$. Further if $G \notin \cE$ and
$\dtd(G) = 4n/7$, then $G \in \cH \cup \cS$.
 \label{clawfree1}
\end{thm}

A proof of Theorem~\ref{clawfree1} is given in
Section~\ref{S:clawfree1}. Since the graphs in the family $\cE$ all have order at most~$10$ and the graphs in the family~$\cS$ all have order at most~$14$, our
second main result, namely Theorem~\ref{clawfree2}, follows as an immediate
consequence of Theorem~\ref{clawfree1}.


\begin{thm}
If $G$ is a connected claw-free graph of order~$n > 14$, then
$\dtd(G) \le 4n/7$, with equality if and only if $G \in \cH$.
 \label{clawfree2}
\end{thm}

\subsection{Preliminary Results and Observations}

Since adding edges to a graph cannot increase its disjunctive total
domination number, we have the following observation.

\begin{ob}
If $H$ is a spanning subgraph of a graph $G$, then $\dtd(G) \le
\dtd(H)$. \label{spanning}
\end{ob}

The following lemma will prove to be useful.

\begin{lem}
If $v$ is a support vertex in a graph $G$ with exactly one neighbor
$w$ that is not a leaf, then there is a $\dtd(G)$-set that
contains~$v$. Further if $d_G(w) = 2$, then there is a $\dtd(G)$-set
that contains both~$v$ and~$w$.
 \label{l:support}
\end{lem}
\textbf{Proof.}  Let $S$ be a $\dtd(G)$-set. If $v \notin S$, then in
order to disjunctively totally dominate the leaf neighbors of $v$, at
least two leaf neighbors of $v$ belong to $S$. But we can simply
replace one of the leaf neighbors of $v$ in $S$ with the vertex $v$.
Hence the set $S$ can be chosen to contain~$v$. Further if $d_G(w) =
2$ and $w \notin S$, then at least one leaf neighbor of $v$ belongs
to $S$ in order to totally dominate or disjunctively totally
dominate~$v$. But we can simply replace such a leaf neighbor of $v$
in $S$ with the vertex~$w$.~\qed

\medskip
The total domination number of a path $P_n$ or a cycle $C_n$ on $n
\ge 3$ vertices is easy to compute (see, \cite{He00}): For $n \ge 3$,
$\gamma_t(C_n)= \lfloor n/2 \rfloor + \lceil n/4 \rceil - \lfloor n/4
\rfloor$. However if we disjunctively total dominate the vertices of
a path or cycle we can do much better.
%
We first establish the value of $\dtd(C_n)$ for a cycle $C_n$.  

\begin{prop}
\label{p:cycle}
For $n \ge 3$, $\dtd(C_{n}) = 2n/5$ if $n \equiv 0 \, (\mod \, 5)$ and $\dtd(C_{n}) = \lceil 2(n+1)/5\rceil$ otherwise.
\end{prop}
\proof We proceed by induction on $n \ge 3$. The result is straightforward to verify for $n \le 9$. Suppose that $n \ge 10$ and that the result is true for all cycles of order less than~$n$. Let $G$ be a cycle on $n$ vertices given by $v_1v_2 \ldots v_nv_1$. We first establish upper bounds on $\dtd(G)$. Let
\[
S = \bigcup_{i = 0}^{ \lfloor n/5 \rfloor - 1} \{v_{5i+1},v_{5i+2}\}.
\]

If $n \equiv 0 \, (\mod \, 5)$, let $D = S$.
If $n \equiv 1 \, (\mod \, 5)$, let $D = S \cup \{v_n\}$.
If $n \equiv i \, (\mod \, 5)$ where $i \in \{2,3,4\}$, let $D = S \cup \{v_{n-4},v_{n-3}\}$.
In all cases, the set $D$ is a DTD-set of $G$. Further if $n \equiv 0 \, (\mod \, 5)$, then $|D| = 2n/5$, while if $n \not\equiv 0 \, (\mod \, 5)$, then $|D| = \lceil 2(n+1)/5 \rceil$. Hence, $\dtd(C_n) \le 2n/5$ if $n \equiv 0 \, (\mod \, 5)$ and $\dtd(C_n) \le \lceil 2(n+1)/5 \rceil$ if $n \not\equiv 0 \, (\mod \, 5)$.

To prove the reverse inequality, let $T$ be a $\dtd(G)$-set. If no two vertices in $T$ are adjacent, then $n$ is even and either $T$ contains all vertices with even subscripts or all vertices with odd subscripts. Hence, $\dtd(G) = |T| = n/2$. Since $n \ge 10$, this contradicts our earlier upper bounds established on $\dtd(G)$. Hence the set $T$ contains two consecutive vertices on the cycle. Renaming vertices if necessary, we may assume that $\{v_1,v_2\} \subseteq T$. We show next that we can choose $T$ so that $T \cap \{v_3,v_4,v_5\} = \emptyset$.

Suppose that $v_3 \in T$. If $v_4 \in T$, let $i$ be the smallest integer such that $i \ge 5$ and $v_i \notin T$, and replace $v_3$ in $T$ with the vertex $v_{i}$. Suppose that $v_4 \notin T$. If $v_5 \in T$, replace $v_3$ in $T$ with $v_4$. If $v_5 \notin T$ and $v_6 \in T$, replace $v_3$ in $T$ with $v_5$. If $v_5 \notin T$ and $v_6 \notin T$, then $v_7 \in T$ and replace $v_3$ in $T$ with $v_6$. In all the above cases, we can choose $T$ so that $v_3 \notin T$.

Suppose that $v_4 \in T$. If $v_5 \notin T$, then $v_6 \in T$ and we can replace $v_4$ in $T$ with the vertex $v_5$. If $v_5 \in T$, let $i$ be the smallest integer such that $i \ge 6$ and $v_i \notin T$, and replace $v_4$ in $T$ with the vertex $v_{i}$. In both cases we can choose $T$ so that $v_4 \notin T$.

Suppose that $v_5 \in T$. Then, $v_6 \in T$. Let $i$ be the smallest integer such that $i \ge 7$ and $v_i \notin T$, and replace $v_5$ in $T$ with the vertex $v_{i}$. Hence we can choose $T$ so that $v_5 \notin T$. Therefore, $T \cap \{v_3,v_4,v_5\} = \emptyset$, implying that $\{v_6,v_7\} \subset T$. Let $T'  = T \setminus \{v_1,v_2\}$, and note that $|T'| = |T| - 2$.

We now let $G'$ be obtained from $G$ by deleting the vertices $v_i$, $1 \le i \le 5$, and adding the edge $v_6v_n$. Then, $G' = C_{n'}$, where $n' = n - 5 \ge 5$. Since $T$ is a DTD-set of $G$, the set $T'$ is a DTD-set of $G'$. Hence, $\dtd(G') \le |T'| = |T| - 2$. Applying the inductive hypothesis to $G'$, we have that $\dtd(G') = 2n'/5 = 2n/5 - 2$ if $n \equiv 0 \, (\mod \, 5)$ and $\dtd(G') = \lceil 2(n'+1)/5 \rceil = \lceil 2(n+1)/5 \rceil - 2$ if $n \not\equiv 0 \, (\mod \, 5)$. This implies that $\dtd(G) = |T| \ge 2n/5$ if $n \equiv 0 \, (\mod \, 5)$ and $\dtd(G) = |T| \ge \lceil 2(n+1)/5 \rceil$ if $n \not\equiv 0 \, (\mod \, 5)$. The desired bounds now follow as a consequence of the upper bounds on $\dtd(G)$ established earlier.~\qed

\medskip
We next establish the value of $\dtd(P_n)$ for a path $P_n$. The proof is similar to that of Proposition~\ref{p:cycle} and is omitted.

\begin{prop}
\label{p:path}
For $n \ge 3$, $\dtd(P_n)) = \lceil 2(n+1)/5\rceil + 1$ if $n \equiv 1 \, (\mod \, 5)$ and $\dtd(P_n) = \lceil 2(n+1)/5\rceil$ otherwise. \end{prop}

We remark that if $G$, for example, is a cycle $C_n$ where $n \equiv
0 \, (\mod \, 15)$, then $\gamma(G) = n/3$ while $\dtd(G) = 2n/5$,
implying that $\dtd(G) - \gamma(G) = n/15$, which can be made
arbitrarily large for sufficiently large~$n$. On the other hand, if
$G$ is obtained from $K_{2,k+2}$, where $k \ge 1$, by adding a
pendant edge to each vertex of degree~$2$ and adding an edge joining
the two vertices of degree~$k+2$, then $\gamma(G) = k+2$ while
$\dtd(G) = 2$, implying that $\gamma(G) - \dtd(G) = k$, which can be
made arbitrarily large for sufficiently large~$k$. Hence we have the
following result.

\begin{ob}
There is no relationship between the domination number $\gamma$ and
the disjunctive total domination $\dtd$ of a graph in the sense that
for every positive integer~$k$, there exist graphs $G$ and $H$ such
that $\dtd(G) - \gamma(G) \ge k$ and $\gamma(H) - \dtd(H) \le k$.
 \label{relate}
\end{ob}

\subsection{Proof of Theorem~\ref{char:graph}}
\label{S:char:graph}

In order to prove Theorem~\ref{char:graph}, we first establish an
upper bound on the disjunctive total domination number of a tree in
terms of its order and we characterize the extremal trees.

\begin{thm}
Let $T \notin \{P_5,P_6\}$ be a tree of order~$n \ge 4$. Then,
$\dtd(T) \le 2(n-1)/3$, with equality if and only if $T \in \cT \cup
\cF \cup \{K_{1,3},T^*\}$.
 \label{char:tree}
\end{thm}
\proof We proceed by induction on the order~$n \ge 4$ of a tree $T
\notin \{P_5,P_6\}$ to show that $\dtd(T) \le 2(n-1)/3$. Furthermore,
if equality holds, then $T \in \cT \cup \cF \cup \{K_{1,3},T^*\}$. If
$n = 4$, then $T = P_4 \in \cT$ or $T = K_{1,3}$. In both cases,
$\dtd(T) = 2 = 2(n-1)/3$. If $n = 5$, then $\dtd(T) = 2 < 2(n-1)/3$
(recall that $T \ne P_5$, by assumption). If $n = 6$, then $\dtd(T)
\le 3 < 2(n-1)/3$ (recall that $T \ne P_6$, by assumption). Hence we
may assume that $n \ge 7$.

If $T$ is a star or a double star, then $\dtd(T) = 2 < 2(n-1)/3$.
Hence we may assume that $\diam(T) \ge 4$.
Let $P$ be a longest path in $T$ and let $P$ be an $r$-$u$ path.
Necessarily, both $r$ and $u$ are leaves. We now root the tree at the
vertex $r$. Let $v$ be the parent of $u$, and let $w$ be the parent
of $v$ and $x$ the parent of $w$ in the rooted tree. We proceed
further with the following series of claims. Recall that $n \ge 7$.

\newpage
\begin{unnumbered}{Claim~I.}
If $d_T(v) \ge 3$, then $\dtd(T) \le 2(n-1)/3$, with equality if and
only if $T = T^* \in \cT$.
\end{unnumbered}
\textbf{Proof.}  Suppose that $d_T(v) \ge 3$. Let $T' = T - u$. Let
$T'$ have order~$n'$, and so $n' = n - 1 \ge 6$. We note that $v$ is
a support vertex in $T'$ with exactly one neighbor that is not a
leaf. By Lemma~\ref{l:support}, there is a $\dtd(T')$-set $S'$ that
contains~$v$. If $T' = P_6$, then $T = T^*$, and we are done. Hence
we may assume that $T' \ne P_6$. Applying the inductive hypothesis to
$T'$, we have that $\dtd(T') \le 2(n'-1)/3 = 2(n-2)/3$, with equality
if and only if $T' \in \cT \cup \cF \cup \{T^*\}$. Since $S'$ is a
DTD-set in $T$, we have that $\dtd(T) \le |S'| \le 2(n-2)/3 <
2(n-1)/3$.~\smallqed

\medskip
By Claim~I, we may assume that $d_T(v) = 2$, for otherwise the
desired result follows. Analogously we may assume that every
(support) vertex at distance~$\diam(T) - 1$ from $r$ in $T$ has
degree~$2$ in $T$.

\begin{unnumbered}{Claim~II.}
If $d_T(w) \ge 3$, then $\dtd(T) \le 2(n-1)/3$, with equality if and
only if $T = F_2 \in \cF$.
\end{unnumbered}
\textbf{Proof.}  Suppose that $d_T(w) \ge 3$. Let $T' = T - \{u,v\}$
have order~$n'$. Then, $n' = n - 2 \ge 5$. If $T' \in \{P_5,P_6\}$,
then either $T = F_2 \in \cF$ or $\dtd(T) < 2(n-1)/3$, and we are
done. Hence we may assume that $T' \notin \{P_5,P_6\}$. Applying the
inductive hypothesis to $T'$, we have that $\dtd(T') \le 2(n'-1)/3 =
2(n-3)/3$, with equality if and only if $T' \in \cT \cup \cF \cup
\{T^*\}$. We proceed further with the following subclaim.

\begin{unnumbered}{Claim~II.1.}
There is a $\dtd(T')$-set that contains $w$ or two neighbors of $w$
in $T'$.
\end{unnumbered}
\textbf{Proof.}
If $w$ has a leaf-neighbor $w'$ in $T'$, then every $\dtd(T')$-set
contains $w$ or two neighbors of $w$ in $T'$ different from~$w'$ in
order to disjunctively totally dominate the leaf $w'$. Hence we may
assume that $w$ has no leaf-neighbor in $T'$, for otherwise the
desired result follows. Then, $d_{T'}(w) \ge 2$ and every neighbor of
$w$ different from its parent $x$ is a support vertex of degree~$2$
in $T'$. Let $v'$ be a child of $w$ in $T'$ and let $u'$ be the child
of $v'$. By Lemma~\ref{l:support}, there is a $\dtd(T')$-set $S'$
that contains~$v'$. If $u' \in S'$, then we simply replace $u'$ in
$S'$ by $w$. Hence in order to disjunctively totally dominate the
vertex $v'$, we have that $w \in S'$ or $S'$ contains two neighbors
of $w$ in $T'$ different from~$v'$. Thus the set $S'$ can be chosen
to contain $w$ or two neighbors of $w$ in $T'$, as desired.~\smallqed

\medskip
By Claim~II.1, there is a $\dtd(T')$-set that contains $w$ or two
neighbors of $w$ in $T'$. Let $S'$ be such a set. Then the set $S'
\cup \{v\}$ is a DTD-set in $T$, and so $\dtd(T) \le |S'| + 1 \le
2(n-3)/3 + 1 < 2(n-1)/3$. This completes the proof of
Claim~II.~\smallqed

\medskip
By Claim~II, we may assume that $d_T(w) = 2$, for otherwise the
desired result follows.  Analogously we may assume that every vertex
at distance~$\diam(T) - 2$ from $r$ in $T$ has degree~$2$ in $T$. Let
$T' = T - \{u,v,w\}$ have order~$n'$. Then, $n' = n - 3 \ge 4$.

Suppose $n' = 4$, and so $n = 7$. If $T' = K_{1,3}$, then either $x$
is the central vertex of $T'$, in which case $\dtd(T) = 3 <
2(n-1)/3$, or $x$ is a leaf of $T'$, in which case $T = T^*$. If $T'
= P_4$, then either $x$ is a central vertex of $T'$, in which case $T
= F_2 \in \cF$, or $x$ is a leaf of $T'$, in which case $T = P_7 \in
\cT$. Hence we may assume that $n' \ge 5$, for otherwise the desired
result follows. If $T' = P_5$, then $n = 8$ and $\dtd(T) = 4 <
2(n-1)/3$. If $T' = P_6$, then $n = 9$ and $\dtd(T) \le 5 <
2(n-1)/3$. Hence we may assume that $T' \notin \{P_5,P_6\}$ (and
still $n' \ge 5$).
Applying the inductive hypothesis to $T'$, we have that $\dtd(T') \le
2(n'-1)/3 = 2(n-4)/3$, with equality if and only if $T' \in \cT \cup
\cF \cup \{T^*\}$. Every $\dtd(T')$-set can be extended to a DTD-set
of $T$ by adding to it the vertices $v$ and $w$, implying that
\begin{equation}
\dtd(T) \le \dtd(T') + 2 \le \frac{2}{3}(n'-1) + 2 = \frac{2}{3}(n-1).
 \label{Eq1}
\end{equation}

This establishes the desired upper bound. Suppose that $\dtd(T) =
2(n-1)/3$. Then we must have equality throughout the Inequality
Chain~(\ref{Eq1}). In particular, $\dtd(T') = 2(n'-1)/3$, implying
that $T' \in \cT \cup \cF \cup \{T^*\}$ (recall that $n' \ge 5$). If
$T' = T^*$, then $n = 10$ and irrespective of which vertex of $T'$
corresponds to the vertex $x$, we have that $\dtd(T) \le 5 <
2(n-1)/3$, a contradiction. Hence, $T' \in \cT \cup \cF$.

Suppose $T' \in \cT$. Then, $T' = T_k$ for some $k \ge 2$, and so $n
= 3k+4$. If $x$ is a leaf of $T'$ or if $x$ is a support vertex of
$T'$, then $\dtd(T) = 2k + 1 < 2(n-1)/3$, a contradiction. If $x$ is
the central vertex of $T'$, then $T = T_{k+1} \in \cT$. If $x$ is
adjacent to the central vertex of $T'$, then $T = F_{k+1} \in \cF$.
Hence we may assume that $T' \in \cF$, for otherwise if $T' \in \cT$,
then $T \in \cT \cup \cF$, as desired. Thus, $T' = F_k$ for some $k
\ge 2$. It is now a routine exercise to check that if $x$ is not the
central vertex of $T'$, then $\dtd(T) = 2k + 1 < 2(n-1)/3$, a
contradiction. Hence, $x$ is the central vertex of $T'$, implying
that $T = F_{k+1} \in \cF$.

Conversely, suppose that $T \in \cT \cup \cF \cup \{K_{1,3},T^*\}$.
If $T \in \{K_{1,3},P_4,T^*\}$, then it is a simple exercise to check
that $\dtd(T) = 2(n-1)/3$. Suppose, therefore, that $T \in \cT \cup
\cF$ has order~$n \ge 7$. Let $S$ be a $\dtd(T)$-set. If $T = T_k \in
\cT$ or if $T = F_k \in \cF$, then by Lemma~\ref{l:support} and its
proof the set $S$ can be chosen to contain all $k$ support vertices
of $T$ and their $k$ neighbors that are not leaf-neighbors, and so
$\dtd(T) = |S| \ge 2(n-1)/3$. As shown earlier, $\dtd(T) \le
2(n-1)/3$ for all trees $T$ of order~$n \ge 4$ different from $P_5$
and $P_6$. Consequently, $\dtd(T) = 2(n-1)/3$. This completes the
proof of Theorem~\ref{char:tree}.~\qed

\medskip
As an immediate consequence of Theorem~\ref{char:tree}, we have the
following result.

\begin{cor}
If $T$ is a tree of order~$n \ge 8$, then $\dtd(T) \le 2(n-1)/3$,
with equality if and only if $T \in \cT \cup \cF$.
 \label{cor:tree}
\end{cor}

We are now in a position to prove Theorem~\ref{char:graph}. Recall
its statement.

\noindent \textbf{Theorem~\ref{char:graph}} \emph{Let $G$ be a
connected graph of order~$n \ge 8$. Then, $\dtd(G) \le 2(n-1)/3$,
with equality if and only if $G \in \cT \cup \cF \cup \cG$.}

\noindent \textbf{Proof.} Let $G $ be a connected graph of order~$n
\ge 8$. If $G$ is a tree, then by Corollary~\ref{cor:tree} we have
that $G \in \cT \cup \cF$, as desired. Hence we may assume that $G$
is not a tree. Among all spanning trees of $G$, let $T$ be chosen so
that $\dtd(T)$ is a minimum. By Observation~\ref{spanning}, $\dtd(G)
\le \dtd(T)$. By Corollary~\ref{cor:tree}, $\dtd(T) \le 2(n-1)/3$.
Consequently, $\dtd(G) \le 2(n-1)/3$.

Suppose that $\dtd(G) = 2(n-1)/3$. Then, $\dtd(T) = 2(n-1)/3$. By our
choice of $T$ and by Corollary~\ref{cor:tree}, this implies that
every spanning tree $T'$ of $G$ satisfies $\dtd(T') = 2(n-1)/3$ and
therefore belongs to the family~$\cT \cup \cF$. In particular, we
note that $n \ge 10$ and $n \equiv 1 \, (\mod \, 3)$. Let $e \in E(G)
\setminus E(T)$.

Suppose that $T \in \cT$. Then, $T = T_k$ for some integer $k \ge 3$.
If $e$ joins a leaf of $T$ to the central vertex of $T$, then
deleting the cycle edge in $T + e$ that is not incident to either end
of $e$ (and therefore joins two vertices of degree~$2$ in $T$) produces a
spanning tree of $G$ that does not belong to the family~$\cT \cup
\cF$, a contradiction. Hence we may assume that $e$ does not join a
leaf of $T$ to the central vertex of $T$. If $e$ does not join two
neighbors of the central vertex in $T$, then deleting a cycle edge of
$T + e$ different from $e$ produces a spanning tree of $G$ that does
not belong to the family~$\cT \cup \cF$, a contradiction. Thus the
edge $e$ joins two neighbors of the central vertex in $T$, and so $T
+ e = G_k \in \cG$. This is true for every edge in $E(G) \setminus
E(T)$. Suppose that there is an edge $f \in E(G) \setminus (E(T) \cup
\{e\})$. Let $e'$ be a cycle edge in $T + e$ different from $e$, and
let $f'$ be a cycle edge in $T + f$ different from $e$ and $e'$.
Then, the tree $(T + \{e,f\}) - \{e',f'\}$ obtained from $T +
\{e,f\}$ by deleting the two edges $e'$ and $f'$ produces a spanning
tree of $G$ that does not belong to the family~$\cT \cup \cF$, a
contradiction. Therefore, $e$ is the only edge of $G$ not in $T$,
implying that $G = T + e = G_k \in \cG$, as desired.

Suppose that $T \in \cF$. Then, $T = F_k$ for some integer $k \ge 3$.
Let $v$ denote the vertex of $T$ that is at distance at least~$3$
from every leaf of $T$ and let $u$ denote the vertex at distance~$2$
from~$v$ in $T$ that is not a support vertex in $T$. If $e$ does not join $u$ and $v$, then deleting a cycle edge in
$T + e$ different from $e$ produces a spanning tree of $G$ that does
not belong to the family~$\cT \cup \cF$, a contradiction. Hence the
edge $e$ must join $u$ and $v$, implying that $E(G) \setminus E(T) =
\{e\}$ and $G = T + e = G_k \in \cG$, as desired.~\qed

\section{Claw-Free Graphs}
\label{S:clawfree1}

In this section, we determine a tight upper bound on the disjunctive
total domination number of a claw-free graph. For this purpose, we
first present some preliminary results and observations. Recall that
the family $\cL = \{L_1,L_2,\ldots,L_{12}\}$ consists of the 12
connected claw-free graphs on seven vertices with total domination
number~$4$. Of these 12 graphs,  six have disjunctive total
domination number~$4$, namely the graphs in the family $\cS_1 =
\{L_1,L_2,L_3,L_5,L_6,L_{10}\}$. Hence since $\dtd(G) \le \gt(G)$ for
all graphs $G$ with no isolated vertex, we have the following
observation.

\begin{ob}
If $G$ is a connected claw-free graph of order~$7$, then $\dtd(G)
\le 4$ with equality if and only if $G \in \cS_1$.
 \label{ob:n7II}
\end{ob}

As a consequence of Corollary~\ref{ob:clawfree}, we have the
following result.

\begin{ob}
If $G$ is a connected claw-free graph of order~$n$ with $\delta(G)
\ge 2$, then $\dtd(G) < 4n/7$, unless $G \in \{C_3,C_7\}$.
\label{ob:clawfree1}
\end{ob}

We remark that every graph of order~$n \ge 8$ in the family $\cT \cup
\cF \cup \cG$ contains a claw, except for the graph $G_3$. Hence we
have the following immediate consequence of Theorem~\ref{char:graph}.

\begin{cor}
If $G$ is a connected claw-free graph of order~$n \ge 8$, then
$\dtd(G) \le 2n/3 - 1$, unless $G = G_3$.
 \label{cor:claw}
\end{cor}

Recall that the set of exceptional graphs is the set $\cE = \{
P_2,P_3,P_5,P_6,C_3,G_3 \}$. We shall need the following lemma.

\begin{lem}
Let $G$ be a connected claw-free graph of order~$n$, where $2 \le n
\le 11$. Then, $\dtd(G) \le 4n/7$, unless $G \in \cE$. Further if $G
\notin \cE$ and $\dtd(G) = 4n/7$, then $G \in \cS_1$.
 \label{l:small}
\end{lem}
\proof If $n = 2$, then $G = P_2$, while if $n = 3$, then $G \in
\{P_3,C_3\}$. If $n = 4$, then $\gt(G) = 2 < 4n/7$. If $n = 5$, then
either $G = P_5$ or $\gt(G) = 2 < 4n/7$. If $n = 6$, then either $G =
P_6$ or $\gt(G) \le 3 < 4n/7$. Since $\{P_2,P_3,P_5,P_6,C_3\} \subset
\cE$ and since $\dtd(G) \le \gt(G)$, it therefore follows that for $2
\le n \le 6$ either $G \in \cE$ or $\dtd(G) < 4n/7$.
If $n = 7$, then by Observation~\ref{ob:n7II}, we have $\dtd(G) \le 4$
with equality if and only if $G \in \cS_1$.
Suppose $8 \le n \le 11$. If $G = G_3$, then $G$ is an
exceptional graph in the family $\cE$. Hence we may assume that $G
\ne G_3$. By Corollary~\ref{cor:claw}, $\dtd(G) \le \lfloor 2n/3 \rfloor - 1$. Since $\lfloor 2n/3 \rfloor - 1 < 4n/7$ for all values of $n \in \{8,9,10,11\}$, we therefore have that $\dtd(G) < 4n/7$.~\qed

\medskip
We are now in a position to present a proof of
Theorem~\ref{clawfree1}. Recall its statement.

\noindent \textbf{Theorem~\ref{clawfree1}} \emph{Let $G$ be a
connected claw-free graph of order~$n \ge 2$. Then, $\dtd(G) \le
4n/7$, unless $G \in \cE$. Further if $G \notin \cE$ and $\dtd(G) =
4n/7$, then $G \in \cH \cup \cS$. }

\noindent \textbf{Proof.} The proof is by induction on $n \ge 2$. If
$n \le 11$, then the desired result follows from Lemma~\ref{l:small}.
This establishes the base cases. Let $n \ge 12$ and assume the result
holds for all connected claw-free graph of order~$n'$, where $2 \le
n' < n$, and let $G$ be a connected claw-free graph of order~$n$.

If $\delta(G) \ge 2$, then by Observation~\ref{ob:clawfree1} (and recalling that $n \ge 12$), we have $\gt(G) < 4n/7$. Hence we may assume that $G$ has a vertex of degree~$1$, for otherwise there is nothing left to prove.
Let $y$ be a leaf and $x$ its neighbor. Let $X = N[x] \setminus
\{y\}$. By the claw-freeness of~$G$, the set $X$ is a clique. We call
a component of $G-X$ a fragment. A fragment isomorphic to a graph
$F$, we call an $F$-fragment. For each fragment $F$ we choose a
vertex $x_F \in X$ that is adjacent to a vertex of $F$. If $w$ is a chosen vertex in $X$ associated with a fragment $F$, then $w = x_F$ and we denote the fragment $F$ by $F_w$.
By the claw-freeness of $G$, every vertex in $X$ is adjacent to vertices
from at most one fragment. The chosen vertices, $x_F$, associated
with a fragment are therefore distinct.

Let $X_1$ be those vertices in $X$ that are not one of the chosen
vertices, $x_F$, where $F \in \cE$. Thus if $v \in X_1$, then either
$v$ is not a chosen vertex associated with a fragment or $v = v_F$
for some fragment $F$ where $F \notin \cE$. Let $Y$ be the set
consisting of the vertices of the $P_1$-fragments together with set
$X_1$. In particular, we note that $\{x,y\} \subseteq Y$, and so $|Y|
\ge 2$. Further if $|Y| \ge 3$, then $|X_1| \ge 2$.

We proceed further with the following algorithm. The algorithm carefully selects a set $S$ of vertices of $G$ depending on the structure of the fragments in $G - X$. The resulting set $S$ is chosen in such a way that it is either a DTD-set of $G$ or can be extended to a DTD-set of $G$. In both cases, the resulting TDT-set of $G$ has cardinality at most four-sevenths the order of $G$.

\begin{unnumbered}{Algorithm~A.}
Let $S$ be the set of vertices in $G$ constructed as follows.
Initially, set $S = \emptyset$. Let $F$ be a fragment of $G$. \\
 1. If $|Y| \ge 4$, then add to $S$ the vertex $x$ and one other
vertex of $X_1$. \\
2. If $|Y| = 2$ or if $|Y| = 3$, then add to $S$ the vertex $x$. \\
3. If $F \ne P_1$ and $F \notin \cE$, then add to $S$ a
$\dtd(F)$-set. \\
4. If $F=P_2$, do the following. \\
\hspace*{0.3cm} 4.1. If both vertices of $F$ are adjacent to $x_F$,
then add $x_F$ to $S$. \\
\hspace*{0.3cm} 4.2. If exactly one vertex of $F$ is adjacent to
$x_F$, then add this vertex of $F$ to $S$. \\
5. If $F = P_3$, do the following. \\
\hspace*{0.3cm} 5.1. If $x_F$ is adjacent to the central vertex of
$F$, then add the vertex $x_F$ and a neighbor \\
\hspace*{1.1cm} of $x_F$ in $F$ to $S$. \\
\hspace*{0.3cm} 5.2. If $x_F$ is not adjacent to the central vertex, then by the claw-freeness of $G$ the vertex \\ \hspace*{1.1cm}  $x_F$ is adjacent to exactly one leaf of $F$. Add the leaf of $F$ adjacent to $x_F$ and the \\ \hspace*{1.1cm}
central vertex of $F$ to $S$. \\
6. If $F = C_3$, then add to $S$ the vertex $x_F$ and a neighbor of
$x_F$ in $F$. \\
7. If $F = P_5$, let $v_1v_2\ldots v_5$ denote the path $F$. \\
\hspace*{0.3cm} 7.1. If $x_F$ is adjacent to a leaf of $F$, say to
$v_1$, then add $x_F$, $v_3$ and $v_4$ to $S$. \\
\hspace*{0.3cm} 7.2. If $x_F$ is not adjacent to a leaf of $F$, then
by the claw-freeness of $G$ the vertex $x_F$ is \\
\hspace*{1.1cm} adjacent to $v_3$ and to at least one of $v_2$ and
$v_4$, say to $v_2$ by symmetry. In this case, \\
\hspace*{1.1cm} add $x_F$, $v_3$ and $v_4$ to $S$. \\
8. If $F = P_6$, let $v_1v_2\ldots v_6$ denote the path $F$. \\
\hspace*{0.3cm} 8.1. If $x_F$ is adjacent to a leaf of $F$, say to
$v_1$, then add $x_F$, $v_4$ and $v_5$ to $S$. \\
\hspace*{0.3cm} 8.2. If $x_F$ is not adjacent to a leaf of $F$ but
adjacent to a support vertex of $F$, say to $v_2$, \\
\hspace*{1.1cm} then by the claw-freeness of $G$ the vertex
$x_F$ is also adjacent to $v_3$. In this case, add \\
\hspace*{1.1cm} $v_2$, $v_4$ and $v_5$ to $S$. \\
\hspace*{0.3cm} 8.3. If $x_F$ is adjacent to neither a leaf nor a
support vertex of $F$, then by the claw-freeness \\
\hspace*{1.1cm} of $G$ the vertex $x_F$ is adjacent to both $v_3$ and
$v_4$. In this case, add $x_F$, $v_2$, $v_4$ and $v_5$ \\
\hspace*{1.1cm} to $S$. \\
9. If $F = G_3$, let $V(F) =
\{u_1,u_2,u_3,v_1,v_2,v_3,w,w_1,w_2,w_3\}$, where $w$ is the central
vertex \\ \hspace*{0.5cm} of $F$, and for $i \in \{1,2,3\}$,
$\{u_i,v_i,w_i\}$ is the set of three vertices at distance~$i$ from
$w$ in \\ \hspace*{0.5cm} $F$, where $u_1u_2u_3$, $v_1v_2v_3$ and
$w_1w_2w_3$ are paths and where $u_1$ and $v_1$ are adjacent. Let \\
\hspace*{0.5cm} $D = \{x_F,u_1,u_2,v_1,v_2,w_1,w_2\}$. \\
\hspace*{0.7cm} 9.1. If $x_F$ is adjacent to $w_3$, then add $(D
\setminus \{w_1,w_2\}) \cup \{w_3\}$ to $S$. \\
\hspace*{0.7cm} 9.2. If $x_F$ is not adjacent to $w_3$ but to $u_3$
or $v_3$, say to $u_3$, then add $(D \setminus \{u_1,u_2\}) \cup \{u_3\}$ \\
\hspace*{1.5cm} to $S$. \\
\hspace*{0.7cm} 9.3. If $x_F$ is adjacent to no leaf of $F$ but is
adjacent to $w_2$, then add $D \setminus \{w_2\}$ to $S$. \\
\hspace*{0.7cm} 9.4. If $x_F$ is adjacent to neither a leaf of $F$
nor to $w_2$ but is adjacent to $u_2$ or $v_2$, say \\
\hspace*{1.5cm} to $u_2$, then add $D \setminus \{u_2\}$ to $S$. \\
\hspace*{0.7cm} 9.5. Suppose $x_F$ is adjacent to neither a leaf nor
a support vertex of $F$. By the \\
\hspace*{1.5cm} claw-freeness of $G$ the vertex $x_F$ is adjacent
to~$w$. \\
\hspace*{1.7cm} 9.5.1. If $x_F$ is adjacent to $w_1$, then add $(D
\setminus \{u_1,w_1\}) \cup \{w\}$ to $S$. \\
\hspace*{1.7cm} 9.5.2. If $x_F$ is not adjacent to $w_1$, then by
the claw-freeness of $G$ the vertex $x_F$ \\
\hspace*{2.9cm}  is adjacent to both $u_1$ and $v_1$, and we add $(D
\setminus \{u_1,v_1\}) \cup \{w\}$ to $S$.

\end{unnumbered}

We make a few comments about Algorithm~A. In Step~3 of the algorithm,
the $\dtd(F)$-set added to $S$ has cardinality at most $4|V(F)|/7$ by
the inductive hypothesis. Further if it has cardinality
exactly~$4|V(F)|/7$, then $F \in \cH \cup \cS$. In Step~8.3 we note
that $G[V(F) \cup \{x_F\}] = H_1$. In Step~9, we note that the set added to $S$ is
a $6$-element subset of $V(F) \cup \{x_F\}$. In Step~3 to Step~9, we
note that the added set contains the vertex $x_F$, except in
Steps~4.2, 5.2, and 8.2. We note that $x \in S$. Let $F_y$ be the
$P_1$-fragment that consists of the vertex~$y$. We proceed further with a series of
claims.

\begin{unnumbered}{Claim~A}
If $|S \cap X| \ge 2$, then $\dtd(G) < 4n/7$.
\end{unnumbered}
\proof Suppose that the set $S$ contains at least two vertices of
$X$. Then by construction the set $S$ is a DTD-set of $G$.
Furthermore in this case the set $S$ is guaranteed to have
cardinality strictly less than~$4n/7$.~\qed

\medskip
By Claim~A, we may assume that $S \cap X = \{x\}$, for otherwise
$\dtd(G) < 4n/7$ and we are done. This implies that the graph $G$ has
the following properties.

\begin{unnumbered}{Claim~B}
Let $F$ be a fragment in $G$. Then the following holds. \\
\indent {\rm (a)} $|Y| \le 3$. \\
\indent {\rm (b)} There is at most one fragment $F$ such that $F \ne
P_1$ and $F \notin \cE$. \\
\indent {\rm (c)} $F \in \{P_1,P_2,P_3,P_6\}$ or $F \notin \cE$. \\
\indent {\rm (d)} If $F$ is a $P_1$-fragment, then $F = F_y$. \\
\indent {\rm (e)} If $F$ is a $P_2$-fragment, then exactly one vertex
of $F$ is adjacent to $x_F$.
\\
\indent {\rm (f)} If $F$ is a $P_3$-fragment, then a leaf, but no
other vertex, of $F$ is adjacent to $x_F$. \\
\indent {\rm (g)} If $F$ is a $P_6$-fragment, then $x_F$ is adjacent
to a support vertex, but not a leaf, of $F$.
\end{unnumbered}

\begin{unnumbered}{Claim~C}
If there is a $P_6$-fragment, then $\dtd(G) < 4n/7$.
\end{unnumbered}
\proof Suppose there is a $P_6$-fragment, $F$. By Claim~B(g), $x_F$
is adjacent to a support vertex, but not a leaf, of $F$. Let
$v_1v_2\ldots v_6$ be the path $F$. We may assume, renaming vertices
if necessary, that the vertex $x_F$ is adjacent to $v_2$ and $v_3$.
Recall that in Step~8.2 of Algorithm~A, we added $v_2$, $v_4$ and
$v_5$ to $S$. We now let $S^* = S \cup \{x_F\}$. By construction the
resulting set $S^*$ is a DTD-set of $G$ since note that $\{x,x_F\}
\subset S^*$. For each fragment $L$ in $G$ different from $F_y$, we
let $V_L = V(L) \cup \{x_L\}$ and note that $|S^* \cap V_L| \le
4|V_L|/7$. Further we note that $|Y| \in \{2,3\}$ but $S^* \cap Y =
\{x\}$, implying that $|S^* \cap Y| < 4|Y|/7$. Hence the DTD-set
$S^*$ is guaranteed to have cardinality strictly less
than~$4n/7$.~\qed

\medskip
By Claim~C, we may assume that there is no $P_6$-fragment, for
otherwise $\dtd(G) < 4n/7$ and we are done. Hence by Claim~B, if $F$
is a fragment different from $F_y$, then $F$ is a $P_2$-fragment or a
$P_3$-fragment or $F \notin \cE$ and $F \ne P_1$. For $i
\in\{1,2,3\}$, let $k_i$ denote the number of $P_i$-fragments in $G$.
By Claim~B(c), $k_1 = 1$.

\begin{unnumbered}{Claim~D}
If every fragment in $G$ different from $F_y$ belongs to~$\cE$, then
$\dtd(G) < 4n/7$.
\end{unnumbered}
\proof Suppose that every fragment in $G$ different from $F_y$
belongs to~$\cE$. Let $F$ be a fragment different from $F_y$. Then,
$F \in \{P_2,P_3\}$. Recall that $k_1 = 1$ and that $n \ge 12$,
implying that $k_2 + k_3 \ge 3$. If $k_2 = 0$, then $k_3 \ge 3$ and
the set $S$ is a DTD-set of $G$. Further, $|S| = 2k_3 + 1$ and $n = 4k_3
+ |Y| \ge 4k_3 + 2$, implying that $\dtd(G) \le |S| \le n/2 < 4n/7$.
Hence we may assume that $k_2 \ge 1$. In this case, we select one
$P_2$-fragment $F$ of $G$ and consider the set $S^* = S \cup
\{x_F\}$. Since $\{x,x_F\} \subset S^*$, we note that $S^*$ is a
DTD-set of $G$. Further, $|S^*| = k_2 + 2k_3 + 2$. Since $n = 3k_2 +
4k_3 + |Y| \ge 3k_2 + 4k_3 + 2$, we have that $4n/7 \ge k_2 + 2k_3 +
2 + 5k_2/7 + 2k_3/7 - 6/7$. Since $k_2 \ge 2$ or $k_2 = 1$ and $k_3
\ge 2$, we note that $5k_2/7 + 2k_3/7 - 6/7 > 0$, implying that $4n/7
> k_2 + 2k_3 + 2 = |S^*| \ge \dtd(G)$. Hence if every fragment in
$G$ different from $F_y$ belongs to~$\cE$, then $\dtd(G) <
4n/7$.~\qed

\medskip
By Claim~D, we may assume that there is a fragment $F$ in $G$ such
that $F \ne F_y$ and $F \notin \cE$. For notational convenience, let $w = x_F$ and let $F_w$ denote this fragment $F$. Since $\{w,x,y\} \subseteq Y$,
we have by Claim~B(a) that $Y = \{w,x,y\}$. Further
by Claim~B(b), $F_w$ is the only fragment different from $F_y$
that does not belong to~$\cE$. By the inductive hypothesis, $\dtd(F_w) \le
4|V(F_w)|/7$ with equality if and only if $F_w \in \cH \cup \cS$. All remaining fragments, if any, different from $F_w$ and $F_y$ are $P_2$- or $P_3$-fragments.

\begin{unnumbered}{Claim~E}
If $k_2 \ge 1$, then $\dtd(G) < 4n/7$.
\end{unnumbered}
\proof Suppose that $k_2 \ge 1$. In this case, we select one
$P_2$-fragment $F$ of $G$ and consider the set $S^* = S \cup
\{x_F\}$. Since $\{x,x_F\} \subset S^*$, we note that $S^*$ is a
DTD-set of $G$. Further, $|S^*| = k_2 + 2k_3 + 2 + \dtd(F_w) \le k_2 + 2k_3 + 2 + 4|V(F_w)|/7$. Since $n = 3k_2 + 4k_3 + 3 + |V(F_w)|$, we have that $4n/7 \ge k_2 + 2k_3 + 2 + 4|V(F_w)|/7 + (5k_2 + 2k_3 - 2)/7 \ge |S^*| + (5k_2 + 2k_3 - 2)/7$. Since $k_2 \ge 1$ and $k_3 \ge 0$, we note that $5k_2 + 2k_3 - 2 > 0$, implying that $\dtd(G) \le |S^*| < 4n/7$.~\qed

\medskip
By Claim~E, we may assume that $k_2 = 0$.

\begin{unnumbered}{Claim~F}
If $k_3 \ge 2$, then $\dtd(G) < 4n/7$.
\end{unnumbered}
\proof Suppose that $k_3 \ge 2$. In this case, the set $S$ is a DTD-set of $G$. Further, $|S| = 2k_3 + 1 + \dtd(F_w) \le 2k_3 + 1 + 4|V(F_w)|/7$. Since $n = 4k_3 + 3 + |V(F_w)|$, we have that $4n/7 \ge 2k_3 + 1 + 4|V(F_w)|/7 + (2k_3 + 5)/7 > 2k_3 + 1 + 4|V(F_w)|/7 \ge |S| \ge \dtd(G)$. Hence, $\dtd(G) \le |S| < 4n/7$.~\qed

\medskip
By Claim~F, we may assume that $k_3 \le 1$.

\begin{unnumbered}{Claim~G}
If $k_3 = 1$, then $\dtd(G) < 4n/7$.
\end{unnumbered}
\proof Suppose that $k_3 = 1$ and consider the $P_3$-fragment $F$ of $G$. For notational convenience, let $z = x_F$ and let $F_z$ denote this fragment $F$. Then, $G[V(F_z) \cup \{z\}]$ is a path $P_4$. Let $z_1z_2z_3z_4$ denote this path, where $z = z_1$. We now consider the graph $G^* = G - \{z_1,z_2,z_3,z_4\}$. We note that $V(G^*) = V(F_w) \cup \{w,x,y\}$. Further since $F_w \notin \cE$, we note that $|V(F_w)| \ge 4$ and therefore $|V(G^*)| \ge 7$. Applying the inductive hypothesis to $G^*$, we have that either $G^* = G_3$ or $\dtd(G^*) \le 4|V(G^*)|/7$.

Suppose $G^* = G_3$. Let $V(G^*) = \{a_1,a_2,a_3,b_1,b_2,b_3,c_1,c_2,c_3,d\}$, where $d$ is the central vertex of $G^*$, and for $i \in \{1,2,3\}$, $\{a_i,b_i,c_i\}$ is the set of three vertices at distance~$i$ from $d$ in $G^*$. Further, $a_1a_2a_3$, $b_1b_2b_3$ and $c_1c_2c_3$ are paths and $a_1$ and $b_1$ are adjacent. In the graph $G^*$, we note that $yxw$ is a path where $d(y) = 1$ and $d(x) = 2$. Thus renaming vertices if necessary, we may assume that either $y = a_3$ (in which case $x = a_2$ and $w = a_1$) or $y = c_3$ (in which case $x = c_2$ and $w = c_1$). If $y = a_3$, then let $S^* = \{a_2,b_1,b_2,c_1,c_2,z_2,z_3\}$. If $y = c_3$, then let $S^* = \{a_1,a_2,b_2,c_2,d,z_2,z_3\}$. In both cases, $S^*$ is a DTD-set of $G$ and $|S^*| = 7$. Thus in this case, we have that $n = 14$ and $\dtd(G) \le |S^*| = n/2 < 4n/7$. Hence we may assume that $G^* \ne G_3$, for otherwise $\dtd(G) < 4n/7$ as desired.

Since $G^* \ne G_3$, we therefore have that $G^* \notin \cE$, and therefore by the inductive hypothesis that $\dtd(G^*) \le 4|V(G^*)|/7$. Every $\dtd(G^*)$ can be extended to a DTD-set of $G$ by adding to it the vertices $z_2$ and $z_3$, implying that $\dtd(G) \le \dtd(G^*) + 2 \le 4|V(G^*)|/7 + 2 = 4(n - 4)/7 + 2 < 4n/7$.~\qed

\medskip
By Claim~G, we may assume that $k_3 = 0$, for otherwise $\dtd(G) <
4n/7$ and we are done. This implies that $d_G(x) = 2$ and $N_G(x) =
\{y,w\}$. Further, $V(G) = V(F_w) \cup \{w,x,y\}$. Since $n \ge 12$,
we note that $|V(F_w)| \ge 9$. Since $F_w \notin \cE$, we also note
that $F_w \ne G_3$. We may assume in what follows that every support
vertex has degree~$2$, for otherwise we can choose the support vertex
$x$ to have degree at least~$3$ and by the above arguments we would
be done. The formal statement follows.

\begin{unnumbered}{Claim~H}
Every support vertex in $G$ has degree~$2$.
\end{unnumbered}

To simplify the notation in what follows, we now rename the vertex
$y$ by $z$, the vertex $x$ by $y$, and the vertex $w$ by $x$. Hence,
$zyx$ is a path in $G$ (corresponding to the previous path $yxw$),
where $z$ is a leaf of $G$ with neighbor $y$ of degree~$2$ in $G$. As
before we define $X = N[x] \setminus \{y\}$ and note by the
claw-freeness of~$G$ that the set $X$ is a clique.

We now consider the components of $G_X = G - (X \cup \{y,z\})$.
Following our earlier notation, we call each such component a
fragment of $G_X$. A fragment isomorphic to a graph $F$, we again
call an $F$-fragment. As before for each fragment $F$ we choose a
vertex $x_F \in X$ that is adjacent to a vertex of $F$. If $w$ is a chosen vertex in $X$ associated with a fragment $F$, then we denote the fragment $F$ by $F_w$.

Suppose that $F$ is a $P_1$-fragment of $G_X$, where $V(F) = \{u\}$. Then, $N(u) \subseteq X \setminus \{x\}$. Suppose $d_G(u) = 1$. Then, $x_F$ is a support vertex with leaf-neighbor $u$. By Claim~H, the vertex $x_F$ has degree~$2$ in $G$. But this implies that $u$ and $x$ are the two neighbors of $x_F$, and therefore that $|X| = 2$. This in turn implies that $d_G(x) = 2$ and that $G$ is a path $P_5$, a contradiction. Hence, $|N(u)| = d_G(u) \ge 2$. By the claw-freeness of $G$, every vertex in $X$ is adjacent to vertices from at most one fragment. Hence we can uniquely associate the set $N(u)$ with the fragment $F$.

As before, let $X_1$ be those vertices in $X$ that are not one of the chosen vertices, $x_F$, where $F \in \cE$. Thus if $v \in X_1$, then either $v$ is not a chosen vertex associated with a fragment or $v = x_F$ for some fragment $F$, where $F \notin \cE$. Let $Y$ be the set consisting of the vertices of the $P_1$-fragments together with the set $X_1 \cup \{x,y,z\}$. In particular, we note that $|Y| \ge 3$.

We now apply Algorithm~A to the graph $G_X$, except we modify the algorithm slightly as follows. In Step~1 and Step~2, we add to $S$ the vertices $x$ and $y$. All other steps in the algorithm remain unchanged, except that we add one additional step, namely Step~10 which states that if $F = P_1$, then add $x_F$ to $S$. We call the resulting modified algorithm, Algorithm~B.

\begin{unnumbered}{Claim~I}
If $|Y| \ge 4$, then $\dtd(G) < 4n/7$ or $G \in \cH$ or $G = L_{13} \in \cS$.
\end{unnumbered}
\proof Suppose that $|Y| \ge 4$. Then, $|X_1| \ge 1$. By construction the set $S$ has cardinality strictly less than~$4n/7$. Further the set $S$ is a DTD-set of $G$ unless $S \cap X = \{x\}$ and there is a $P_2$-fragment $F$ with exactly one vertex adjacent to $x_F$. Hence, we may assume that $S \cap X = \{x\}$ and that there is a $P_2$-fragment $F$ with exactly one vertex adjacent to $x_F$, for otherwise $\dtd(G) \le |S| < 4n/7$. In particular, we note that there is no fragment $F$ with $x_F \in S$. Thus, if $F$ is a fragment in $G_X$ and $F \in \cE$, then $F$ is either a $P_2$-fragment or a $P_3$-fragment or a $P_6$-fragment, while if $F$ is a fragment in $G_X$ and $F \notin \cE$, then $F \ne P_1$, and so $|V(F)| \ge 4$. Further if $F$ is a $P_2$-fragment, then exactly one vertex of $F$ is adjacent to $x_F$; if $F$ is a $P_3$-fragment, then a leaf, but no other vertex, of $F$ is adjacent to $x_F$; if $F$ is a $P_6$-fragment, then $x_F$ is adjacent to a support vertex, but not a leaf, of $F$.

We now consider a $P_2$-fragment $F$. As observed earlier, exactly one vertex of $F$ is adjacent to $x_F$. Let $x_F = x_1$ and let $F_1 = F$. Further let $V(F_1) =
\{y_1,z_1\}$ where $x_1y_1z_1$ is a path. Suppose $z_1$ is not a leaf in $G$. Then, $z_1$ is adjacent to a vertex $x_2 \in X$, where $x_2 \ne x_1$. By the claw-freeness of $G$, we can uniquely associate $\{x_1,x_2\}$ with the $P_2$-fragment $F$. In this case, replacing the vertex $y_1$ in $S$ with the two vertices $x_1$ and $x_2$ produces a DTD-set of $G$ of cardinality strictly less than~$4n/7$. Hence we may assume that $z_1$ is a leaf in $G$, for otherwise $\dtd(G) < 4n/7$. By Claim~H, the support vertex~$y_1$ therefore has degree~$2$ in $G$. Let $w$ be an arbitrary vertex in $X_1$ and let $T = S \cup \{w\}$. Then, $T$ is a DTD-set of $G$.

Let $\ell_1$, $\ell_2$ and $\ell_3$ denote the number of $P_2$-, $P_3$-, and $P_6$-fragments in $G$, respectively. Let $F_X$ be the disjoint union of all fragments in $G_X$ that do not belong to $\cE$. Applying the inductive hypothesis to each component of $F_X$, we have that $\dtd(F_X) \le 4|V(F_X)|/7$. Thus, $|T| = 3 + \ell_1 + 2\ell_2 + 3\ell_3 + \dtd(F_X) \le 3 + \ell_1 + 2\ell_2 + 3\ell_3 + 4|V(F_X)|/7$. Since $\ell_1 \ge 1$, $\ell_2 \ge 0$ and $\ell_3 \ge 0$, we note that $(5\ell_1 + 2\ell_2 - 5)/7 + \ell_3 \ge 0$, with equality if and only if $\ell_1 = 1$ and $\ell_2 = \ell_3 = 0$. Recall that $|Y| \ge 4$.

Hence,
\[
\begin{array}{lcl} \2
\displaystyle{4n/7 } & = & \displaystyle{ 4( 3\ell_1 + 4\ell_2 + 7\ell_3 + |Y| + |V(F_X)| )/7 } \\ \2
& \ge & \displaystyle{ 4( 3\ell_1 + 4\ell_2 + 7\ell_3 + 4 + |V(F_X)| )/7 }  \\ \2
& = & \displaystyle{ (3 + \ell_1 + 2\ell_2 + 3\ell_3 + 4|V(F_X)|/7) + (5\ell_1 + 2\ell_2 - 5)/7 + \ell_3  } \\ \2
& \ge & \displaystyle{ |T| + (5\ell_1 + 2\ell_2 - 5)/7 + \ell_3  } \\ \2
& \ge & |T| \\ \2
& \ge & \dtd(G).
\end{array}
\]

If $\dtd(G) < 4n/7$, then there is noting left to prove. Hence we may assume that $\dtd(G) = 4n/7$. This implies that we have equality throughout the above inequality chain. Hence, $|Y| = 4$, and so $Y = \{w,x,y,z\}$ and $X_1 = \{w\}$. Further, $\ell_1 = 1$, $\ell_2 = \ell_3 = 0$, $\dtd(F_X) = 4|V(F_X)|/7$, and $\dtd(G) = |T|$. Thus, $F_1$ is the only fragment in~$\cE$, implying that $d_G(x) = 3$ and $N_G(x) = \{x_1,y,w\}$. Since $|X_1| = 1$, there is at most one fragment in $G_X$ not in~$\cE$. However since $n \ge 12$, there is at least one fragment in $G_X$ not in~$\cE$. Consequently, there is exactly one fragment $F$ in $G_X$ not in~$\cE$ and $x_F = w$. We note that $F = F_X$ and $\dtd(F) = 4|V(F)|/7$.
Let $H = G - V(F_w)$ and note that $H = H_1 \in \cH$ and $V(H) = \{w,x,x_1,y,y_1,z,z_1\}$, where $w$ is the central vertex of $H$ with neighbors $x$ and $x_1$ that are adjacent and where $xyz$ and $x_1y_1z_1$ are paths emanating from $x$ and $x_1$, respectively. Recall that $\{w,x,y,y_1\} \subset T$.

As observed earlier, $\dtd(F) = 4|V(F)|/7$. Applying the inductive hypothesis to $F$, we have that $F \in \cH \cup \cS$. We proceed further with the following two subclaims. Recall that by our earlier assumption, $\dtd(G) = |T| = 4n/7$. Further recall that $T = S \cup \{w\}$ and that $\{x,y,y_1\} \subset T$.

\begin{unnumbered}{Claim~I.1}
If $F \in \cH$, then $G \in \cH$.
\end{unnumbered}
\proof Suppose $F \in \cH$. Then, $F = H_t$ for some $t \ge 1$. Let $v$ be an arbitrary vertex of the complete graph $K_t$ used to construct $F$ and let $H_v$ denote the copy of $H_1$ that contains~$v$. Let $V(H_v) = \{u_1,u_2,u_3,v,v_1,v_2,v_3\}$, where $vu_1v_1$ is the triangle in $H_1$ and where $u_1u_2u_3$ and $v_1v_2v_3$ are paths in $H_v$. Renaming vertices if necessary, we may assume that the vertex $w$ is adjacent to at least one vertex in $V(H_v)$ in $G$. Recall that $S \cap V(F)$ was chosen to be an arbitrary $\dtd(F)$-set. We may assume that $S \cap V(F)$ consists of the support vertices of $F$ and the neighbors of the support vertices that are not leaves. In particular, we note that $\{u_1,u_2,v_1,v_2\} \subseteq S$.

Suppose that $w$ is adjacent to $u_3$ or $v_3$, say to $v_3$ by symmetry. Then, $T \setminus \{v_1,v_2\}$ is a DTD-set of $G$, contradicting the fact that $\dtd(G) = |T|$. Hence, $w$ is not adjacent to a leaf of $H_v$ in $G$. By Claim~H, this implies that both $u_2$ and $v_2$ have degree~$2$ in $G$ and are therefore not adjacent to~$w$. Suppose that $w$ is adjacent to $u_1$ or $v_1$, say to $v_1$ by symmetry. By the claw-freeness of $G$, we have that $w$ is adjacent to $v$ and to $u_1$. But then $T \setminus \{v_1\}$ is a DTD-set of $G$, contradicting the minimality of $T$. Hence, the vertex $v$ is the only vertex of $H_v$ adjacent to~$w$. By the claw-freeness of $G$, this implies that $w$ is adjacent to every vertex of the complete graph $K_t$ used to construct $F$. But then $G = H_{t+1} \in \cH$. This completes the proof of Claim~I.1.~\smallqed

\begin{unnumbered}{Claim~I.2}
If $F \in \cS$, then $G = L_{13} \in \cS$ or $G = H_2 \in \cH$.
\end{unnumbered}
\proof Suppose that $F \in \cS$, where we recall that $\cS = \{L_1,L_2,L_3,L_5,L_6,L_{10},L_{13},L_{14}\}$. We consider each of the eight possibilities for $F$ in turn. Let $A = \{w,x,y,y_1\}$ and let $B = S \cap V(F)$. Then, $B$ is a $\dtd(F)$-set and $|B|= 4|V(F)|/7$. Further, $T = A \cup B$. We note that if $F \in \cS_1 = \{L_1,L_2,L_3,L_5,L_6,L_{10}\}$, then $|V(F)| = 7$ and $|B|=4$.

Suppose that $F = L_1$. Let $F$ be the path $v_1v_2 \ldots v_7$. Suppose that $w$ is adjacent to $v_1$ or $v_7$, say to $v_1$ by symmetry. Then, $(T \setminus B) \cup \{v_4,v_5,v_6\}$ is a DTD-set of $G$, contradicting the fact that $\dtd(G) = |T|$. Hence, $w$ is not adjacent to a leaf of $F$ in $G$. Thus both leaves of $F$ are leaves in $G$. By Claim~H, this implies that both $v_2$ and $v_6$ have degree~$2$ in $G$ and are therefore not adjacent to~$w$. By the claw-freeness of $G$, the vertex $w$ is adjacent to at least one of $v_3$ and $v_5$. By symmetry, we may assume that $w$ is adjacent to $v_3$. By the claw-freeness of $G$, we have that $w$ is adjacent to $v_4$. If $w$ is also adjacent to $v_5$, then $(T \setminus B) \cup \{v_2,v_4,v_6\}$ is a DTD-set of $G$, contradicting the fact that $\dtd(G) = |T|$. Hence, in this case $w$ is not adjacent to $v_5$. But then the graph $G$ is determined and $G = L_{13} \in \cS$.

Suppose that $F = L_2$. Let $V(F)=\{v, a_1, a_2, a_3, b_1, b_2, b_3\}$ where $v$ is the central vertex of $F$ and for $i \in\{1,2,3\}$, $\{a_i, b_i\}$ denotes the pair of vertices at distance $i$ from $v$ in $F$, with $a_1a_2a_3$ and $b_1b_2b_3$ as paths and $a_1$ adjacent to $b_1$. Suppose that $w$ is adjacent to $a_3$ or $b_3$, say to $a_3$ by symmetry. Then, $(T \setminus B) \cup \{v,b_1,b_2\}$ is a DTD-set of $G$, contradicting the fact that $\dtd(G) = |T|$. Hence, $w$ is not adjacent to a leaf of $F$ in $G$. Thus both leaves of $F$ are leaves in $G$. By Claim~H, this implies that both $a_2$ and $b_2$ have degree~$2$ in $G$ and are therefore not adjacent to~$w$. Suppose that $w$ is adjacent to $a_1$ or $b_1$, say to $a_1$ by symmetry. By the claw-freeness of $G$, we have then that $w$ is adjacent to both $v$ and $b_1$. Then, $(T \setminus B) \cup \{v,a_2,b_2\}$ is a DTD-set of $G$, a contradiction. Thus, $v$ is the only vertex of $F$ adjacent to $w$ in $G$, implying that $G = H_2 \in \cH$.

Suppose that $F = L_3$. Let $V(F) = \{a_1,a_2,a_3, a_4,a_5,a_6,v\}$, where $F$ is obtained from the path $a_1a_2 \ldots a_6$ by adding a new vertex $v$ and adding the edges $va_4$ and $va_5$. If $w$ is adjacent to $a_1$, let $D = (T \setminus B) \cup \{v,a_4\}$. If $w$ is adjacent to both $a_3$ and $a_4$, let $D = (T \setminus B) \cup \{v,a_2,a_4\}$. In all other cases (noting by Claim~H that if $a_1$ is a leaf of $G$, then $a_2$ is not adjacent to $w$), let $D = (T \setminus B) \cup \{a_2,a_3,a_5\}$. Then, $D$ is a DTD-set of $G$ with $|D| < |T|$, contradicting the fact that $\dtd(G) = |T|$.

Suppose that $F = L_5$. Let $F$ be obtained from the path $a_1a_2 \ldots a_7$ by adding the edge $a_5a_7$. If $w$ is adjacent to $a_1$, let $D = (T \setminus B) \cup \{a_4,a_5\}$. If $w$ is adjacent to both $a_3$ and $a_4$, let $D = (T \setminus B) \cup \{a_2,a_4,a_5\}$. In all other cases (noting by Claim~H that if $a_1$ is a leaf of $G$, then $a_2$ is not adjacent to $w$), let $D = (T \setminus B) \cup \{a_2,a_3,a_5\}$. Then, $D$ is a DTD-set of $G$ with $|D| < |T|$, contradicting the fact that $\dtd(G) = |T|$.

Suppose that $F = L_6$. Let $V(F)=\{a_1, a_2, a_3, v_1, v_2, v_3, v_4\}$, where $F$ is obtained from the path $a_1a_2a_3v_1v_2v_3v_4$ by adding the edges $v_1v_4$ and $v_2v_4$. If $w$ is adjacent to $a_1$, let $D = (T \setminus B) \cup \{v_1,v_2\}$. If $w$ is adjacent to both $a_3$ and $v_1$, let $D = (T \setminus B) \cup \{a_2,v_1,v_2\}$. In all other cases, let $D = (T \setminus B) \cup \{a_2,a_3,v_1\}$. Then, $D$ is a DTD-set of $G$ with $|D| < |T|$, contradicting the fact that $\dtd(G) = |T|$.

Suppose that $F = L_{10}$. Let $F$ be the cycle $v_1v_2 \ldots v_7v_1$. By the claw-freeness of $G$, the vertex $w$ is adjacent to two consecutive vertices on the cycle $F$. Renaming vertices, if necessary, we may assume that $w$ is adjacent to $v_1$ and $v_2$. Then, $(T \setminus B) \cup \{v_4,v_5\}$ is a DTD-set of $G$, contradicting the fact that $\dtd(G) = |T|$.

Suppose that $F = L_{13}$. Let $F$ be obtained from the disjoint union of the paths $a_1a_2 \ldots a_7$ and $b_1b_2 \ldots b_6$ by adding a new vertex $v$ and joining $v$ to $a_3$, $a_4$, $b_3$ and $b_4$. Let $S_F = \{a_2, a_3, a_5, a_6, b_2, b_3, b_5, v\}$ and note that $|S_F| = 8 = |B| = 4|V(F)|/7$. If $w$ is adjacent to $a_1$, let $D = (T \setminus B) \cup (S_F \setminus \{a_2,a_3\})$. If $w$ is adjacent to $a_7$, let $D = (T \setminus B) \cup (S_F \setminus \{a_3,a_5,a_6\}) \cup \{a_4\}$. If $w$ is adjacent to $a_5$, let $D = (T \setminus B) \cup (S_F \setminus \{a_5\})$. If $w$ is adjacent to $a_3$, let $D = (T \setminus B) \cup (S_F \setminus \{a_3\})$. If $w$ is adjacent to $b_1$ or $b_6$, say to $b_6$ by symmetry, let $D = (T \setminus B) \cup (S_F \setminus \{b_5\})$. If $w$ is adjacent to $b_3$ or $b_4$, let $D = (T \setminus B) \cup (S_F \setminus \{b_3\})$. By the claw-freeness of $G$ and by Claim~H, this exhausts all possibilities. In all cases, $D$ is a DTD-set of $G$ with $|D| < |T|$, contradicting the fact that $\dtd(G) = |T|$.

Suppose that $F = L_{14}$. Let $F$ be obtained from the disjoint union of the paths $a_1a_2 \ldots a_7$ and $b_1b_2 \ldots b_7$ by forming a clique on the set $\{a_3,a_4,b_3,b_4\}$. Let $S_F = \{a_2, a_3, a_5, a_6, b_2, b_3$, $b_5, b_6\}$ and note that $|S_F| = 8 = |B| = 4|V(F)|/7$. If $w$ is adjacent to $a_7$ or $b_7$, say to $a_7$ by symmetry, let $D = (T \setminus B) \cup (S_F \setminus \{a_6\})$. If $w$ is adjacent to $a_1$ or $b_1$, say to $a_1$ by symmetry, let $D = (T \setminus B) \cup (S_F \setminus \{a_2\})$. If $w$ is adjacent to $a_5$ or $b_5$, say to $a_5$ by symmetry, let $D = (T \setminus B) \cup (S_F \setminus \{a_3,a_5\}) \cup \{a_4\}$. If $w$ is adjacent to $a_3$ or $b_3$, say to $a_3$ by symmetry, let $D = (T \setminus B) \cup (S_F \setminus \{a_3\})$. By the claw-freeness of $G$ and by Claim~H, this exhausts all possibilities. In all cases, $D$ is a DTD-set of $G$ with $|D| < |T|$, contradicting the fact that $\dtd(G) = |T|$. This completes the proof of Claim~I.2.~\smallqed

\medskip
Claim~I now follows from Claim~I.1 and Claim~I.2.~\qed

\medskip
By Claim~I, we may assume that $|Y| = 3$, for otherwise the desired result follows. Hence, $Y = \{x,y,z\}$, implying that every vertex in $X \setminus \{x\}$ is associated with a fragment $F$ where $F \in \cE$. Recall that the set $S$ is obtained by Algorithm~B (as explained in the paragraph immediately preceding Claim~I) and $\{x,y\} \subset S$. Recall that $n \ge 12$. If $|S \cap X| \ge 3$, then the set $S \setminus \{x\}$ is a DTD-set of $G$ and by construction, $|S \setminus \{x\} | < 4n/7$. Hence we may assume that $|S \cap X| \le 2$.

\begin{unnumbered}{Claim~J}
If $|S \cap X| = 2$, then $\dtd(G) < 4n/7$ or $G \in \{L_{13},L_{14}\} \subset \cS$.
\end{unnumbered}
\proof Suppose that $|S \cap X| = 2$. Then the set $S$ is a DTD-set of $G$. Let $x'$ be the vertex in $S \cap X$ different from $x$ and let $F$ be the fragment associated with the vertex $x'$, and so $x' = x_F$. Let $S_F = S \cap V(F)$. Since $x_F$ is the only vertex of $X \setminus \{x\}$ that belongs to $S$, we note that if $L$ is a fragment in $G_X$ different from $F$, then either $L$ is a $P_2$-fragment with exactly one vertex of $L$ adjacent to $x_L$ or $L$ is a $P_3$-fragment with a leaf, but no other vertex, of $L$ adjacent to $x_L$ or $L$ is a $P_6$-fragment and $x_L$ is adjacent to a support vertex, but not a leaf, of $L$. However by Claim~H, every support vertex in $G$ has degree~$2$, implying that there is no $P_6$-fragment.

Suppose that $F = G_3$ and let $F$ be described as in Step~9 of Algorithm~A. If $|X| \ge 3$, then $|S| < 4n/7$. Therefore we may assume in this case that $X = \{x,x_F\}$, and so $V(G) = \{x,y,z\} \cup V(F)$ and $n = 14$. If $S_F$ is constructed in Step~9.1 of Algorithm~A, let $T = (S \setminus \{x_F,w_3\}) \cup \{w_1\}$. If $S_F$ is constructed in Step~9.2, let $T = (S \setminus \{u_3,v_1,x_F\}) \cup \{u_1,w\}$. By Claim~H, we note that $S_F$ cannot be constructed in Step~9.3 or Step~9.4. If $S_F$ is constructed in Step~9.5.1, then $x_F$ is adjacent to both $w$ and $w_1$. By the claw-freeness of $G$, we note that in this case, $x_F$ is adjacent to neither $u_1$ nor $v_1$, and so $G = L_{13} \in \cS$. If $S_F$ is constructed in Step~9.5.2, then $G = L_{14} \in \cS$. In all cases, either $G \in \{L_{13},L_{14}\}$ or $G \notin \{L_{13},L_{14}\}$ and the constructed set $T$ is a DTD-set of $G$ with $|T| \le 7 < 4n/7$. Hence in what follows we may assume that $F \ne G_3$, for otherwise the desired result follows. If there is a fragment $L$, different from $F$, such that $L = G_3$, then by Algorithm~A we note that $x_L \in S$, a contradiction since $S \cap X = \{x,x_F\}$ and $x_L \ne x_F$. Hence there is no fragment isomorphic to $G_3$.

If $S_F$ is not constructed in Step~8.3 of Algorithm~A, then $|S| < 4n/7$. Hence we may assume that $S_F$ is constructed  in Step~8.3 of Algorithm~A. In particular, we note that $S_F = \{v_2,v_4,v_5,x_F\}$ and $|S_F| = 4 = 4|V(F)|/7$. Since $n \ge 12$, there are at least two fragments. Let $L$ be a fragment different from $F$ and let $S_L = S \cap V(L)$. If there are three or more fragments, or if $F$ and $L$ are the only two fragments but $S_L$ is not constructed as in Step~5.2 of Algorithm~A, then $|S| < 4n/7$. Hence we may assume that $F$ and $L$ are the only two fragments and that $L$ is a $P_3$-fragment with a leaf, but no other vertex, of $L$ adjacent to $x_L$. This implies that $G = L_{13} \in \cS$.~\qed

\medskip
Recall that since $|Y| = 3$, every vertex in $X \setminus \{x\}$ is associated with a fragment $F$ where $F \in \cE$. By Claim~J, we may assume that $S \cap X = \{x\}$, and so for every fragment $F$ we have $x_F \notin S$. This implies that if $F$ is a fragment in $G_X$, then either $F$ is a $P_2$-fragment with exactly one vertex of $F$ adjacent to $x_L$ (see Step~4.2 in Algorithm~A) or $F$ is a $P_3$-fragment with a leaf, but no other vertex, of $F$  adjacent to $x_F$ (see Step~5.2 in Algorithm~A) or $F$ is a $P_6$-fragment and $x_F$ is adjacent to a support vertex, but not a leaf, of $F$ (see Step~8.2 in Algorithm~A). However by Claim~H, every support vertex in $G$ has degree~$2$, implying that there is no $P_6$-fragment.

Let $\ell_1$ and $\ell_2$ denote the number of $P_2$- and $P_3$-fragments in $G$, respectively. Since $n \ge 12$, we note that $\ell_1 + \ell_2 \ge 3$. Further, $n = 3\ell_1 + 4\ell_2 + 3$. Let $w$ be an arbitrary vertex from the set $X \setminus \{x\}$ and let $S^* = S \cup \{w\}$. We note that the set $S^*$ is a DTD-set of $G$ and $|S^*| = \ell_1 + 2\ell_2 + 3$.

Suppose $\ell_1 = 0$. Then the set $S$ is a DTD-set of $G$ and $\ell_2 \ge 3$. Further, $|S| = 2\ell_2+2$ and $4n/7 = 4(4\ell_2 + 3)/7 = (2\ell_2 + 2) + (2\ell_2 - 2)/7 \ge |S| + 4/7 > |S| \ge \dtd(G)$. Hence we may assume that $\ell_1 \ge 1$, for otherwise $\dtd(G) < 4n/7$.

Suppose $\ell_1 \ge 2$. Then, $4n/7 = 4(3\ell_1 + 4\ell_2 + 3)/7 = (\ell_1 + 2\ell_2 + 3) + (5\ell_1 + 2\ell_2 - 9)/7 \ge |S^*| + 3/7 > |S^*| \ge \dtd(G)$. Hence we may assume that $\ell_1 = 1$, for otherwise $\dtd(G) < 4n/7$.

Since $\ell_1 = 1$ and $\ell_1 + \ell_2 \ge 3$, we note that $\ell_2 \ge 2$. Suppose $\ell_2 \ge 3$. Then, $4n/7 = 4(3\ell_1 + 4\ell_2 + 3)/7 = (\ell_1 + 2\ell_2 + 3) + (5\ell_1 + 2\ell_2 - 9)/7 \ge |S^*| + 2/7 > |S^*| \ge \dtd(G)$. Hence we may assume that $\ell_2 = 2$, for otherwise $\dtd(G) < 4n/7$.

Since $\ell_1 = 1$ and $\ell_2 = 2$, we have that $G = L_{14} \in \cS$. This completes the proof of Theorem~\ref{clawfree1}.~\qed

\medskip

\end{document}